\documentclass[11pt]{article}
\usepackage{enumerate}
\usepackage{amssymb,a4wide,latexsym,makeidx,epsfig,fleqn}
\usepackage{amsthm}
\usepackage{amsmath}
\usepackage{enumerate}
\usepackage{graphicx}
\usepackage{epstopdf}
\newtheorem{theorem}{Theorem}[section]
\newtheorem{remark}[theorem]{Remark}
\newtheorem{example}[theorem]{Example}
\newtheorem{definition}[theorem]{Definition}
\newtheorem{lemma}[theorem]{Lemma}

\begin{document}
\textwidth 150mm \textheight 225mm
\title{The left row rank of quaternion unit gain graphs in terms of girth\footnote{This work is supported by the National Natural Science Foundations of China (No. 12371348, 12201258), the Postgraduate Research \& Practice Innovation Program of Jiangsu Normal University (No. 2024XKT1702).}}
\author{{ Yong Lu\footnote{Corresponding author.}, Qi Shen, JiaXu Zhong}\\
{\small  School of Mathematics and Statistics, Jiangsu Normal University,}\\ {\small  Xuzhou, Jiangsu 221116,
People's Republic
of China.}\\
{\small E-mails: luyong@jsnu.edu.cn, sq\_jsnu@163.com, JXZhong@163.com}}

\date{}
\maketitle
\begin{center}
\begin{minipage}{120mm}
\vskip 0.3cm
\begin{center}
{\small {\bf Abstract}}
\end{center}
{\small Let $\Phi=(G,U(\mathbb{Q}),\varphi)$  be a quaternion unit gain graph (or $U(\mathbb{Q})$-gain graph). The adjacency matrix of $\Phi$ is denoted by $A(\Phi)$ and the left row rank of $\Phi$ is denoted by $r(\Phi)$.
If $\Phi$ has at least one cycle, then the length of the shortest cycle in $\Phi$ is the girth of $\Phi$, denoted by $g$.
In this paper, we prove that $r(\Phi)\geq g-2$ for $\Phi$. Moreover, we characterize $U(\mathbb{Q})$-gain graphs satisfy $r(\Phi)=g-i$ ($i=0,1,2$) and all quaternion unit gain graphs with rank 2.
The results will generalize the corresponding results of simple graphs (Zhou et al. Linear Algebra Appl. (2021), Duan et al. Linear Algebra Appl. (2024) and Duan, Discrete Math. (2024)), signed graphs (Wu et al. Linear Algebra Appl. (2022)), and complex unit gain graphs (Khan, Linear Algebra Appl. (2024)).

 }
\end{minipage}
\end{center}

\section{Introduction }
Let $G=(V(G),E(G))$ be a simple graph with vertex set $V(G)=\{v_{1},v_{2},\ldots,v_{n}\}$ and edge set $E(G)$.
The \emph{adjacency matrix} $A(G)$ of $G$ is the symmetric $n\times n$ matrix with entries $a_{ij}=1$ if $v_{i}$ is adjacent to $v_{j}$ and $a_{ij}=0$ otherwise.
The \emph{rank} (resp., \emph{nullity}) of $G$ is the rank of $A(G)$, denoted by $r(G)$ (resp., $\eta(G)$).

Let $e_{ij}$ be the \emph{oriented edge} from $v_{i}$ to $v_{j}$, and let $\overrightarrow{E}(G)$ be the set of oriented edges such that $e_{ij},e_{ji}\in\overrightarrow{E}(G)$.
A \emph{gain graph} is a triple $\Phi=(G,\Omega,\varphi)$ consisting of an \emph{underlying graph} $G=(V,E)$, the \emph{gain group} $\Omega$ and the \emph{gain function} $\varphi:\overrightarrow{E}(G)\rightarrow\Omega$ such that $\varphi(e_{ij})=\varphi(e_{ji})^{-1}=\overline{\varphi(e_{ji})}$.
The \emph{adjacency matrix} $A(\Phi)$ of $\Phi$ is the Hermitian $n\times n$ matrix with entries $h_{ij}=\varphi(e_{ij})=\varphi_{v_{i}v_{j}}$ if $e_{ij}\in\overrightarrow{E}(G)$ and $h_{ij}=0$ otherwise.
The \emph{rank} of $\Phi$ is the rank of $A(\Phi)$, denoted by $r(\Phi)$.

Collatz et al. \cite{CS} first proposed to characterize all graphs of order $n$ with $r(G)<n$. Until today, this problem is still unsolved.
This problem is of great interest in both chemistry and mathematics.
For a bipartite graph $G$ which corresponds to an alternant hydrocarbon in chemistry, if $r(G)<n$, it is indicated that the corresponding molecule is unstable.
The rank of a graph is also meaningful in mathematics since it is related to the singularity of adjacency matrix.
In the past decades, researchers have used graph parameters to characterize the boundaries of the ranks (or nullities) of graphs and to identify the extremal graphs.
Recently, researchers have extended the results of simple graphs (when $\varphi(\overrightarrow{E})=\{1\}$) to signed graphs (when $\varphi(\overrightarrow{E})\subseteq\{1,-1\}$), which we can refer to \cite{HHL, LWZ, lwn, WS, wlt},
and they have also extended these results to complex unit gain graphs (when $\varphi(\overrightarrow{E})\subseteq\{T\in \mathbb{C}: |T|=1\}$), which we can refer to \cite{hhd, hhy, landy, LUWH, lwx, LWZ1, landwu, REFF, xzw, yqt}.

Let $\mathbb{Q}$ be a four-dimensional vector space over $\mathbb{R}$ (which is the field of the real numbers) with an ordered basis, denoted by $1$, $i$, $j$ and $k$.
A vector $q=x_{0}+x_{1}i+x_{2}j+x_{3}k\in \mathbb{Q}$ is a \emph{real quaternion} (or \emph{quaternion}),
where $x_{0},x_{1},x_{2},x_{3}$ are real numbers and $i,j,k$ satisfy the following conditions:
$$i^{2}=j^{2}=k^{2}=-1;$$
$$ij=-ji=k, jk=-kj=i, ki=-ik=j.$$
According to \cite{ZF}, if $x$ and $y$ are two different quaternions, then $xy\neq yx$, in general.

Let $q=x_{0}+x_{1}i+x_{2}j+x_{3}k\in \mathbb{Q}$.
The \emph{conjugate} $\bar{q}$ (or $q^{\ast}$) of $q$ is $\bar{q}=x_{0}-x_{1}i-x_{2}j-x_{3}k$.
The \emph{modulus} of $q$ is $|q|=\sqrt{q\bar{q}}=\sqrt{x_{0}^{2}+x_{1}^{2}+x_{2}^{2}+x_{3}^{2}}$.
If $q\neq 0$, then the \emph{inverse} of $q$ is $q^{-1}=\frac{\bar{q}}{|q|^{2}}$.
The \emph{real part} of $q$ is $Re(q)=x_{0}$.
The \emph{imaginary part} of $q$ is $Im(q)=x_{1}i+x_{2}j+x_{3}k$.
The \emph{ left (right) row rank} of a quaternion matrix $A\in \mathbb{Q}^{m\times n}$ is the maximum number of rows of $A$ that are left (right) linearly independent.
The \emph{ left (right) column rank} of a quaternion matrix $A\in \mathbb{Q}^{m\times n}$ is the maximum number of columns of $A$ that are left (right) linearly independent.

If the gain group $\Omega$ is $U(\mathbb{Q})=\{q\in \mathbb{Q}: |q|=1\}$, then $\Phi=(G,U(\mathbb{Q}),\varphi)$ (or $\widetilde{G}$ for short) is a \emph{quaternion unit gain graph} (or $U(\mathbb{Q})$-gain graph).
In this paper, the \emph{rank} $r(\widetilde{G})$ of a quaternion matrix $A(\widetilde{G})\in \mathbb{Q}^{n\times n}$ is defined to be the left row rank of $A(\widetilde{G})$.
Belardo et al. \cite{BBCR} defined the adjacency, Laplacian and incidence matrices for a $U(\mathbb{Q})$-gain graph and studied their properties.
Zhou and Lu \cite{QNZ} obtained the relationship between the left row rank of a $U(\mathbb{Q})$-gain graph and the rank of its underlying graph.
Kyrchei et al. \cite{K} provided a combinatorial description of the determinant of the Laplacian matrix of a $U(\mathbb{Q})$-gain graph by using row-column noncommutative determinants.

For $u,v\in V(\widetilde{G})$, we denote by $d_{\widetilde{G}}(u,v)$ the length of the shortest path between $u$ and $v$ in $\widetilde{G}$.
If $u\in V(\widetilde{G})$ and $\widetilde{H}$ is a subgraph of $\widetilde{G}$, then we denote by $d_{\widetilde{G}}(u,\widetilde{H})$ the length of the shortest path between $u$ and $\widetilde{H}$ in $\widetilde{G}$.
For a vertex $x\in V(\widetilde{G})$, we define $N_{\widetilde{G}}(x)$ as the \emph{neighbor set} of $x$ in $\widetilde{G}$.
The \emph{degree} of $x$ is denoted by $d_{\widetilde{G}}(x)$ and $d_{\widetilde{G}}(x)=|N_{\widetilde{G}}(x)|$.
If $d_{\widetilde{G}}(x)=1$, then $x$ is called a \emph{pendant vertex}.
If $\widetilde{G}$ has at least one cycle, then the length of the shortest cycle in $\widetilde{G}$ is the \emph{girth} of $\widetilde{G}$, denoted by $g$.

We use $\widetilde{P}_{n}$, $\widetilde{C}_{n}$ and $\widetilde{S}_{n}$ to denote a $U(\mathbb{Q})$-gain path, a $U(\mathbb{Q})$-gain cycle and a $U(\mathbb{Q})$-gain star on $n$ vertices, respectively.
By $\widetilde{K}_{n}$, $\widetilde{K}_{a,b}$ and $\widetilde{K}_{r,s,t}$, we
refer to a $U(\mathbb{Q})$-gain complete graph, a $U(\mathbb{Q})$-gain complete bipartite graph and a $U(\mathbb{Q})$-gain complete tripartite graph, respectively.
Let $\widetilde{G}$ be a $U(\mathbb{Q})$-gain graph with vertex set $V(\widetilde{G})$ and $U\subseteq V(\widetilde{G})$, we denote by $\widetilde{G}-U$ the $U(\mathbb{Q})$-gain graph obtained from $\widetilde{G}$ by removing the vertices in $U$ together with all incident edges. When $U=\{x\}$, we write $\widetilde{G}-U$ simply as $\widetilde{G}-x$. Sometimes we use the notation $\widetilde{G}-\widetilde{H}$ instead of $\widetilde{G}-V(\widetilde{H})$ when $\widetilde{H}$ is an induced subgraph of $\widetilde{G}$. If $\widetilde{G}_{1}$ is an induced subgraph of $\widetilde{G}$ and $x$ is a vertex not in $\widetilde{G}_{1}$, then we write the subgraph of $\widetilde{G}$ induced by $V(\widetilde{G}_{1})\cup\{x\}$ simply as $\widetilde{G}_{1}+x$.

Let $\xi: V(\widetilde{G})\rightarrow U(\mathbb{Q})$ be a switching function. Switching the $U(\mathbb{Q})$-gain graph $\widetilde{G}$ by $\xi$ means forming a new quaternion unit gain graph $\widetilde{G}^{\xi}$, whose underlying graph is the same as $\widetilde{G}$, but whose gain function is defined on an edge $uv$ by $\varphi^{\xi}(uv)= \xi(u)^{-1}\varphi(uv)\xi(v)$.
If there exists a switching function $\xi$ such that $\widetilde{G}_{2}=\widetilde{G}_{1}^{\xi}$, then $\widetilde{G}_{1}$ and $\widetilde{G}_{2}$ are called \emph{switching equivalent}, denoted by $\widetilde{G}_{1}\leftrightarrow \widetilde{G}_{2}$.
Note that two switching equivalent $U(\mathbb{Q})$-gain graphs have the same rank.

Let $\widetilde{C}_{p}$ and $\widetilde{C}_{q}$ be two vertex-disjoint cycles with $v\in V(\widetilde{C}_{p})$ and $u\in V(\widetilde{C}_{q})$.
Let $\widetilde{\infty}(p,l,q)$ (as shown in Fig. $1$) be obtained from $\widetilde{C}_{p}$ and $\widetilde{C}_{q}$ by connecting $v$ and $u$ with a path $\widetilde{P}_{l}$.
When $l=1$, the graph $\widetilde{\infty}(p,1,q)$ (as shown in Fig. $1$) is obtained from $\widetilde{C}_{p}$ and $\widetilde{C}_{q}$ by identifying $v$ with $u$.
Let $\widetilde{P}_{p+2},\widetilde{P}_{l+2},\widetilde{P}_{q+2}$ be three paths, where $p,l,q\geq0$ and at most one of $p,l,q$ is 0. Let $\widetilde{\theta}(p,l,q)$ (as shown in Fig. $1$) be obtained from $\widetilde{P}_{p+2}$, $\widetilde{P}_{l+2}$ and $\widetilde{P}_{q+2}$ by identifying the three initial vertices and the three terminal vertices.
Note that any $\widetilde{\infty}(p,l,q)$-graph (resp., $\widetilde{\theta}(p,l,q)$-graph) is obtained from $\widetilde{\infty}(p,l,q)$ (resp., $\widetilde{\theta}(p,l,q)$) by attaching trees to some of its vertices.

\begin{figure}[htbp]
 \centering
 \includegraphics[scale=0.75]{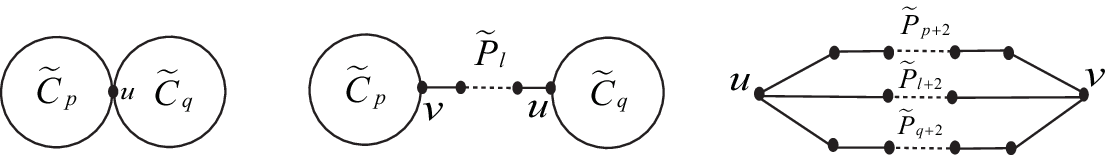}
 \caption{$\widetilde{\infty}(p,1,q)$, $\widetilde{\infty}(p,l,q)$ and $\widetilde{\theta}(p,l,q)$.}
\end{figure}

Here, we introduce some results about girth $g$. For a simple graph, Cheng et al. \cite{CL} obtained the upper bound of $\eta(G)$ in terms of girth $g$,
\begin{align*}
\eta(G)\leq\left\{\begin{array}{ll}
n-g+2,  ~~if~4\mid g;\\
n-g, ~~~~~~~if\;4\nmid g.\\
\end{array}\right.
\end{align*}
Zhou et al. \cite{ZQ} characterized the extremal graphs that satisfy $\eta(G)=n-g$ or $\eta(G)=n-g+2$.
Connected graphs $G$ with $\eta(G)=n-g-1$ have been characterized by Chang et al. \cite{candl}.

Duan et al. \cite{DF} proved $p(G)\geq\lceil\frac{g}{2}\rceil-1$, where $p(G)$ is the number of positive eigenvalues of $G$. The extremal graphs corresponding to $p(G)=\lceil\frac{g}{2}\rceil-1$ and $p(G)=\lceil\frac{g}{2}\rceil$ have been completely characterized, respectively.
Duan \cite{Df} proved $i_-(G)\geq\lceil\frac{g}{2}\rceil-1$, where $i_-(G)$ is the number of negative eigenvalues of $G$, and they characterized the extremal graphs that satisfy $i_-(G)=\lceil\frac{g}{2}\rceil-1$ or $i_-(G)=\lceil\frac{g}{2}\rceil$.
Wu et al. \cite{wlt} extended the results of Zhou et al. (\cite{ZQ}) to signed graphs.
Khan \cite{SK} extended these results to complex unit gain graphs.

In this paper, we extend these results to $U(\mathbb{Q})$-gain graphs.
In Section $2$, we will give some lemmas about $U(\mathbb{Q})$-gain graphs.
In Section $3$, we will prove that $r(\widetilde{G})\geq g-2$, and characterize the extremal $U(\mathbb{Q})$-gain graphs that satisfy $r(\widetilde{G})=g-2$.
In Sections $4$ and $5$, we will characterize the extremal $U(\mathbb{Q})$-gain graphs with $r(\widetilde{G})=g-1$ or $r(\widetilde{G})=g$, and all $U(\mathbb{Q})$-gain graphs with rank 2. The results will generalize the corresponding results of simple graphs (Zhou et al. \cite{ZQ}, Duan et al. \cite{DF} and Duan \cite{Df}), signed graphs (Wu et al. \cite{wlt}), and complex unit gain graphs (Khan \cite{SK}).

\section{Preliminaries}

In this section, we list some known results.

\noindent\begin{lemma}\label{le:2.1}\cite{QNZ}
Let $\widetilde{P}_{n}$ be a $U(\mathbb{Q})$-gain path. If $n$ is odd, then $r(\widetilde{P}_{n})=n-1$; if $n$ is even, then $r(\widetilde{P}_{n})=n$.
\end{lemma}

\noindent\begin{definition}\label{de:2.2}\cite{QNZ}
Let $\widetilde{C}_{n} (n\geq3)$ be a $U(\mathbb{Q})$-gain cycle and let $$\varphi(\widetilde{C}_{n})=\varphi_{v_{1}v_{2}}\varphi_{v_{2}v_{3}}\cdots \varphi_{v_{n-1}v_{n}}\varphi_{v_{n}v_{1}}.$$ Then $\widetilde{C}_{n}$ is said to be:
\begin{displaymath}
\left\{\
        \begin{array}{ll}
          \rm Type~1,&  \emph{if}~\varphi(\widetilde{C}_{n})=(-1)^{n/2}~\emph{and}~n~\emph{is~even};\\
          \rm Type~2,& \emph{if}~\varphi(\widetilde{C}_{n})\neq(-1)^{n/2}~\emph{and}~n~\emph{is~even};\\
          \rm Type~3,& \emph{if}~Re\left((-1)^{{(n-1)}/{2}}\varphi(\widetilde{C}_{n})\right)\neq 0~\emph{and}~n~\emph{is~odd};\\
          \rm Type~4,& \emph{if}~Re\left((-1)^{{(n-1)}/{2}}\varphi(\widetilde{C}_{n})\right)=0~\emph{and}~n~\emph{is~odd}.
        \end{array}
      \right.
\end{displaymath}
\end{definition}

\noindent\begin{lemma}\label{le:2.3}\cite{QNZ}
Let $\widetilde{C}_{n}$ be a $U(\mathbb{Q})$-gain cycle of order $n$. Then
\begin{align*}
r(\widetilde{C}_{n})=\left\{\begin{array}{ll}
n-2, & if \;\widetilde{C}_{n}\; is \;of\; Type \;1;\\
n,& if \;\widetilde{C}_{n}\; is \;of\; Type \;2 \;or\; 3;\\
n-1,& if \;\widetilde{C}_{n}\; is \;of\; Type \;4.\\
\end{array}\right.
\end{align*}
\end{lemma}

\noindent\begin{lemma}\label{le:2.4}\cite{QNZ}
Let $\widetilde{G}$ be a $U(\mathbb{Q})$-gain graph and let $v$ be a vertex of $\widetilde{G}$. Then $r(\widetilde{G})-2\leq r(\widetilde{G}-v)\leq r(\widetilde{G})$.
\end{lemma}

According to the Lemma \ref{le:2.4}, we can obtain the following lemma.

\noindent\begin{lemma}\label{le:2.5}
Let $\widetilde{G}$ be a $U(\mathbb{Q})$-gain graph and let $\widetilde{G}_{1}$ be an induced subgraph of $\widetilde{G}$. Then $r(\widetilde{G}_{1})\leq r(\widetilde{G})$.
\end{lemma}

\noindent\begin{lemma}\label{le:2.6}\cite{QNZ}
Let $\widetilde{G}$ be a $U(\mathbb{Q})$-gain graph and let $x$ be a pendant vertex of $\widetilde{G}$. If $y$ is adjacent to  $x$ in $\widetilde{G}$, then $r(\widetilde{G})=r(\widetilde{G}-x-y)+2$.
\end{lemma}

\noindent\begin{lemma}\label{le:2.7}\cite{QNZ}
Let $\widetilde{G}=\widetilde{G}_{1}\cup \widetilde{G}_{2}\cup \cdots \cup \widetilde{G}_{t}$, where $\widetilde{G}_{1}, \widetilde{G}_{2},\ldots, \widetilde{G}_{t}$ are connected components of a $U(\mathbb{Q})$-gain graph $\widetilde{G}$. Then $r(\widetilde{G})=\sum^{t}_{i=1}r(\widetilde{G}_{i})$.
\end{lemma}

\noindent\begin{lemma}\label{le:2.13}\cite{QNZ}
Let $\widetilde{G}$ be a $U(\mathbb{Q})$-gain graph obtained by identifying a vertex of a $U(\mathbb{Q})$-gain cycle $\widetilde{C}_{n}$ with a vertex of a $U(\mathbb{Q})$-gain graph $\widetilde{G}_{1}$ of order $m(m\geq 1)$ $(V(\widetilde{C}_{n})\cap V(\widetilde{G}_{1})=\{u\})$ and let $\widetilde{G}_{2}=\widetilde{G}_{1}-u$. Then
\begin{align*}
\left\{\begin{array}{ll}
r(\widetilde{G})=n-2+r(\widetilde{G}_{1}), & if \;\widetilde{C}_{n}\; is \;of\; Type \;1;\\
r(\widetilde{G})=n+r(\widetilde{G}_{2}),& if \;\widetilde{C}_{n}\; is \;of\; Type \;2;\\
r(\widetilde{G})=n-1+r(\widetilde{G}_{1}),& if \;\widetilde{C}_{n}\; is \;of\; Type \;4;\\
n-1+r(\widetilde{G}_{2})\leq r(\widetilde{G})\leq n+r(\widetilde{G}_{1}),& if \;\widetilde{C}_{n}\; is \;of\; Type \;3.\\
\end{array}\right.
\end{align*}
\end{lemma}

A \emph{dominating set} for a $U(\mathbb{Q})$-gain graph $\widetilde{G}$ is a subset $D$ of $V(\widetilde{G})$ such that every vertex not in $D$ is adjacent to at least one vertex of $D$.

\noindent\begin{lemma}\label{le:2.8}
Let $\widetilde{G}$ be a connected $U(\mathbb{Q})$-gain graph and let $\widetilde{H}$ be an induced subgraph of $\widetilde{G}$. If $r(\widetilde{G})=r(\widetilde{H})$, then $V(\widetilde{H})$ is a dominating set of $\widetilde{G}$.
\end{lemma}
\noindent\textbf{Proof.}
On the contrary, we find a vertex $v\in V(\widetilde{G})\backslash V(\widetilde{H})$ such that $N_{\widetilde{H}}(v)=\emptyset$. Let $\widetilde{P}_{t}$ ($t\geq3$) be the shortest path from $v$ to $\widetilde{H}$ and let $\widetilde{G}_{1}=\widetilde{H}+\widetilde{P}_{t}$.

If $t$ is even, by Lemmas \ref{le:2.4} and \ref{le:2.6}, then $r(\widetilde{G}_{1})\geq r(\widetilde{H})+t-2>r(\widetilde{G})$, a contradiction.

If $t$ is odd, by Lemma \ref{le:2.6}, then $r(\widetilde{G}_{1})=r(\widetilde{H})+t-1>r(\widetilde{G})$, a contradiction.

So $|N_{\widetilde{H}}(v)|\geq 1$ for any $v\in V(\widetilde{G})\backslash V(\widetilde{H})$. Thus $V(\widetilde{H})$ is a dominating set of $\widetilde{G}$.
\quad $\square$\\

\noindent\begin{lemma}\label{le:2.9}
Let $\widetilde{G}$ be a $U(\mathbb{Q})$-gain graph and let $\widetilde{C}_{g}$ be a shortest cycle in $\widetilde{G}$. If there exists $v\in V(\widetilde{G})\backslash V(\widetilde{C}_{g})$ such that $|N_{\widetilde{C}_{g}}(v)|\geq2$, then $g=3$ or $4$.
\end{lemma}
\noindent\textbf{Proof.}
Let $v_{1}, v_{2}\in N_{\widetilde{C}_{g}}(v)$ and let $\widetilde{P}_{a},\widetilde{P}_{b}$ be two paths in $\widetilde{C}_{g}$ from $v_{1}$ to $v_{2}$ ($a\leq b$). Then $a-1\leq\frac{g}{2}$ and $\widetilde{P}_{a}+v$ is a cycle in $\widetilde{G}$. Thus $(a-1)+2\geq g$, i.e., $g\leq4$. Hence $g=3$ or $4$.\quad $\square$\\

If two pendant vertices have the same neighbor in $\widetilde{G}$, then they are called \emph{pendant twins}.

\noindent\begin{lemma}\label{le:2.10}
Let $v_{1},~v_{2}$ be pendant twins of a $U(\mathbb{Q})$-gain graph $\widetilde{G}$. Then $r(\widetilde{G})=r(\widetilde{G}-v_{1})=r(\widetilde{G}-v_{2})$.
\end{lemma}
\noindent\textbf{Proof.}
Let $V(\widetilde{G})=\{v_{1},v_{2},\ldots,v_{n}\}$ and let $v_{1},v_{2}$ be pendant twins with the same neighbor $v_{3}$ in $\widetilde{G}$. Let $h_{13}=\varphi_{v_{1}v_{3}}$, $h_{23}=\varphi_{v_{2}v_{3}}$.
\begin{align*}
r(\widetilde{G})= r\left (
\begin{array}{ccccccccccccc}
 0 & 0 & h_{13} & \textbf{0}_{1\times(n-3)}\\
 0 & 0 & h_{23} & \textbf{0}_{1\times(n-3)}\\
 h_{31} & h_{32} & 0 & \alpha\\
 \textbf{0}_{(n-3)\times1} & \textbf{0}_{(n-3)\times1} & \overline{\alpha} & \textbf{0}_{(n-3)\times(n-3)}\\
 \end{array}
 \right),
 \end{align*}
where $\alpha=(\varphi_{v_{3}v_{4}}, \varphi_{v_{3}v_{5}},\cdots,\varphi_{v_{3}v_{n}})$.
We multiply $-h_{13}h_{32}$ on the left side of 2-th row and add it to 1-th row,
or multiply $-h_{23}h_{31}$ on the left side of 1-th row and add it to 2-th row.
\begin{align*}
r(\widetilde{G})&= r\left (
\begin{array}{ccccccccccccc}
 0 & 0 & 0 & \textbf{0}_{1\times(n-3)}\\
 0 & 0 & h_{23} & \textbf{0}_{1\times(n-3)}\\
 h_{31} & h_{32} & 0 & \alpha\\
 \textbf{0}_{(n-3)\times1} & \textbf{0}_{(n-3)\times1} & \overline{\alpha} & \textbf{0}_{(n-3)\times(n-3)}\\
 \end{array}
 \right)\\
 &= r\left (
\begin{array}{ccccccccccccc}
 0 & 0 & h_{13} & \textbf{0}_{1\times(n-3)}\\
 0 & 0 & 0 & \textbf{0}_{1\times(n-3)}\\
 h_{31} & h_{32} & 0 & \alpha\\
 \textbf{0}_{(n-3)\times1} & \textbf{0}_{(n-3)\times1} & \overline{\alpha} & \textbf{0}_{(n-3)\times(n-3)}\\
 \end{array}
 \right).
 \end{align*}
Thus $r(\widetilde{G})=r(\widetilde{G}-v_{1})=r(\widetilde{G}-v_{2})$.\quad $\square$\\

From Lemmas \ref{le:2.4}, \ref{le:2.6} and \ref{le:2.10}, we get the following three lemmas.

\noindent\begin{lemma}\label{le:2.10.}
Let $\widetilde{G}_{0}$ be a $U(\mathbb{Q})$-gain graph of order $n-p$.
Let $\widetilde{G}_{1}$ be obtained from $\widetilde{G}_{0}$ and $\widetilde{S}_{p}$ by joining a vertex of $\widetilde{G}_{0}$ with the center of $\widetilde{S}_{p}$.
Let $\widetilde{G}_{2}$ be obtained from $\widetilde{G}_{0}$ and $\widetilde{S}_{p+1}$ by identifying a vertex of $\widetilde{G}_{0}$ with the center of $\widetilde{S}_{p+1}$.
Then $r(\widetilde{G}_{2})\leq r(\widetilde{G}_{1})\leq r(\widetilde{G}_{2})+2$.
\end{lemma}

\noindent\begin{lemma}\label{le:2.11}
Let $\widetilde{G}_{0}$ be a $U(\mathbb{Q})$-gain graph of order $n-p-q$ and let $v_{1},v_{2}\in V(\widetilde{G}_{0})$.
Let $\widetilde{G}_{1}$ be obtained from $\widetilde{G}_{0}$, $\widetilde{S}_{p+1}$ and $\widetilde{S}_{q+1}$ by identifying $v_{1}$ with the center of $\widetilde{S}_{p+1}$, $v_{2}$ with the center of $\widetilde{S}_{q+1}$.
Let $\widetilde{G}_{2}$ be obtained from $\widetilde{G}_{0}$ and $\widetilde{S}_{p+q+1}$ by identifying $v_{1}$ with the center of $\widetilde{S}_{p+q+1}$.
Then $r(\widetilde{G}_{2})\leq r(\widetilde{G}_{1})\leq r(\widetilde{G}_{2})+2$.
\end{lemma}

\noindent\begin{lemma}\label{le:2.12}
Let $\widetilde{G}_{1}$ and $\widetilde{G}_{2}$ be two $U(\mathbb{Q})$-gain graphs and let $v_{1}\in V(\widetilde{G}_{1})$, $v_{2}\in V(\widetilde{G}_{2})$.
Let $\widetilde{G}_{3}$ be obtained from $\widetilde{G}_{1}$ and $\widetilde{G}_{2}$ by connecting $v_{1}$ and $v_{2}$ with a path $\widetilde{P}_{t}$ $(t\geq3)$.
Let $\widetilde{G}_{4}$ be obtained from $\widetilde{G}_{1}$ and $\widetilde{G}_{2}$ by identifying $v_{1}$ with $v_{2}$ and adding $t-1$ pendant vertices to $v_{1}$.
Then $r(\widetilde{G}_{3})\geq r(\widetilde{G}_{4})$.
\end{lemma}

\section{$U(\mathbb{Q})$-gain graphs $\widetilde{G}$ with $r(\widetilde{G})=g-2$}
In this section, we will prove that $r(\widetilde{G})\geq g-2$ for all $U(\mathbb{Q})$-gain graphs and characterize the graphs which satisfy  $r(\widetilde{G})=g-2$. First, we will prove the following lemma.

\noindent\begin{lemma}\label{le:3.1}
Let $\widetilde{K}_{a,b}~(a,b\geq2)$ be a $U(\mathbb{Q})$-gain complete bipartite graph and let $V(\widetilde{K}_{a,b})=V_{1}\cup V_{2}$, $|V_{1}|=a$, $|V_{2}|=b$. Then $r(\widetilde{K}_{a,b})=2$ if and only if all the $\widetilde{C}_{4}$ in $\widetilde{K}_{a,b}$ are of Type 1.
\end{lemma}
\noindent\textbf{Proof.}
\textbf{Necessity:}
Let\begin{align*}
 &A(\widetilde{K}_{a,b})=\left (
 \begin{array}{ccccccc}
 \textbf{0} & A_{1}\\
 A_{1}^{\ast} & \textbf{0}\\
 \end{array}
 \right),
  \end{align*}
where $ A_{1}^{\ast}$ is the conjugate transpose of $A_{1}$.
Let $\alpha_{1},\alpha_{2},\ldots,\alpha_{a}$ be the row vectors of $A_{1}$.
Since $r(\widetilde{K}_{a,b})=2$, we have $r(A_{1})=1$.
Without loss of generality, let $\alpha_{i}=k_{i}\alpha_{1}$ and $k_{i}\neq0~(i=2,3,\ldots,a)$.

Let $u_{1},u_{2}\in V_{1}$ and $v_{1},v_{2}\in V_{2}$.
For convenience, we assume that $\alpha_{1},\alpha_{2}$ are the vectors corresponding to $u_{1},u_{2}$ in $A_{1}$, respectively.
Let $\varphi_{u_{i}v_{j}}$ be the elements in $A_{1}$ corresponding to the edge $u_{i}v_{j}$ $(i,j\in\{1,2\})$.
Since $\alpha_{2}=k_{2}\alpha_{1}$, we have
$$\varphi_{u_{2}v_{1}}=k_{2}\varphi_{u_{1}v_{1}}~\textrm{and}~\varphi_{u_{2}v_{2}}=k_{2}\varphi_{u_{1}v_{2}}.$$
Let $\widetilde{C}_{4}$ be the $U(\mathbb{Q})$-gain cycle with vertices $u_{1},v_{1},u_{2},v_{2}$ in turn.
Then
\begin{align*}
\varphi(\widetilde{C}_{4})&=\varphi_{u_{1}v_{1}}\varphi_{v_{1}u_{2}}\varphi_{u_{2}v_{2}}\varphi_{v_{2}u_{1}}\\
&=\varphi_{u_{1}v_{1}}(k_{2}\varphi_{u_{1}v_{1}})^{-1}(k_{2}\varphi_{u_{1}v_{2}})\varphi_{v_{2}u_{1}}\\
&=\varphi_{u_{1}v_{1}}\varphi_{u_{1}v_{1}}^{-1}k_{2}^{-1}k_{2}\varphi_{u_{1}v_{2}}\varphi_{v_{2}u_{1}}\\
&=(\varphi_{u_{1}v_{1}}\varphi_{u_{1}v_{1}}^{-1})(k_{2}^{-1}k_{2})(\varphi_{u_{1}v_{2}}\varphi_{v_{2}u_{1}})=1.
\end{align*}
By Definition \ref{de:2.2}, $\widetilde{C}_{4}$ is of Type $1$.

\textbf{Sufficiency:} Let $A_{1}$, $\alpha_{1},\alpha_{2},\ldots,\alpha_{a}$ be the same described in the proof of ``Necessity",
$V_{1}=\{u_{1},u_{2},\ldots,u_{a}\}$ and $V_{2}=\{v_{1},v_{2},\ldots,v_{b}\}$.
For convenience, we assume that $\alpha_{i}$ is the vector corresponding to $u_{i}$ in $A_{1}$ $(i=1,2,\ldots,a)$.
The induced subgraph with vertex set $\{u_{1},u_{2},v_{1},v_{2}\}$ is $\widetilde{C}_{4}$.
Since all the $\widetilde{C}_{4}$ in $\widetilde{K}_{a,b}$ are of Type 1, we have $\varphi(\widetilde{C}_{4})=\varphi_{u_{1}v_{1}}\varphi_{v_{1}u_{2}}\varphi_{u_{2}v_{2}}\varphi_{v_{2}u_{1}}=1.$
So
\begin{align*}
\varphi_{v_{1}u_{2}}\varphi_{u_{2}v_{2}}\varphi_{v_{2}u_{1}}&=\varphi_{v_{1}u_{1}}\\
\varphi_{u_{2}v_{2}}\varphi_{v_{2}u_{1}}&=\varphi_{u_{2}v_{1}}\varphi_{v_{1}u_{1}}.
\end{align*}
Let $k_{2}=\varphi_{u_{2}v_{2}}\varphi_{v_{2}u_{1}}=\varphi_{u_{2}v_{1}}\varphi_{v_{1}u_{1}}$, then $\varphi_{u_{2}v_{1}}=k_{2}\varphi_{u_{1}v_{1}}$ and $\varphi_{u_{2}v_{2}}=k_{2}\varphi_{u_{1}v_{2}}$.

The induced subgraph with vertex set $\{u_{1},u_{2},v_{1},v_{3}\}$ is $\widetilde{C}_{4}$.
Since all the $\widetilde{C}_{4}$ in $\widetilde{K}_{a,b}$ are of Type 1, we have $\varphi(\widetilde{C}_{4})=\varphi_{u_{1}v_{1}}\varphi_{v_{1}u_{2}}\varphi_{u_{2}v_{3}}\varphi_{v_{3}u_{1}}=1.$
So
\begin{align*}
\varphi_{v_{1}u_{2}}\varphi_{u_{2}v_{3}}\varphi_{v_{3}u_{1}}&=\varphi_{v_{1}u_{1}}\\
\varphi_{u_{2}v_{3}}\varphi_{v_{3}u_{1}}&=\varphi_{u_{2}v_{1}}\varphi_{v_{1}u_{1}}.
\end{align*}
Then $k_{2}=\varphi_{u_{2}v_{3}}\varphi_{v_{3}u_{1}}$, thus $\varphi_{u_{2}v_{3}}=k_{2}\varphi_{u_{1}v_{3}}.$

Since the induced subgraph with vertex set $\{u_{1},u_{2},v_{1},v_{j}\}$ is $\widetilde{C}_{4}$,
we use the same method to get $\varphi_{u_{2}v_{j}}=k_{2}\varphi_{u_{1}v_{j}}~(j=1,2,\ldots,b)$. Thus $\alpha_{2}=k_{2}\alpha_{1}$.

Repeat the above steps, we can get that $\alpha_{i}=k_{i}\alpha_{1}~(i=2,3,\ldots,a)$, i.e., $r(A_{1})=1$.
Thus $r(\widetilde{K}_{a,b})=2$.
\quad $\square$\\

Now, we will prove one of the main results of this paper.

\noindent\begin{theorem}\label{th:3.2}
Let $\widetilde{G}$ be a connected $U(\mathbb{Q})$-gain graph with girth $g$. Then $r(\widetilde{G})\geq g-2$, where the equality holds if and only if $\widetilde{G}$ is one of the following types:
\begin{enumerate}[(a)]
\item $\widetilde{G}=\widetilde{C}_{g}$ is of Type $1$;
\item $\widetilde{G}=\widetilde{K}_{a,b}$ for $a,b\geq2$ and all the $\widetilde{C}_{4}$ in $\widetilde{K}_{a,b}$ are of Type $1$.
\end{enumerate}
\end{theorem}

\noindent\textbf{Proof.}
Let $\widetilde{G}$ be a connected $U(\mathbb{Q})$-gain graph and let $\widetilde{C}_{g}$ be the shortest cycle in $\widetilde{G}$.
By Lemmas \ref{le:2.3} and \ref{le:2.5}, we have $r(\widetilde{G})\geq r(\widetilde{C}_{g})\geq g-2$. Next, we prove that the equation is true.

``If'' part: If $\widetilde{G}=\widetilde{C}_{g}$ is of Type $1$, by Lemma \ref{le:2.3}, then $r(\widetilde{G})=g-2$.

If $\widetilde{G}=\widetilde{K}_{a,b}$ for $a,b\geq2$ and all the $\widetilde{C}_{4}$ in $\widetilde{K}_{a,b}$ are of Type $1$, by Lemma \ref{le:3.1}, then $g=4$ and $r(\widetilde{G})=2=g-2$.

``Only if'' part: Since $g-2=r(\widetilde{G})\geq r(\widetilde{C}_{g})\geq g-2$, we have $r(\widetilde{C}_{g})=g-2$. By Lemma \ref{le:2.3} and Definition \ref{de:2.2}, $\widetilde{C}_{g}$ is of Type $1$ and $g$ is even.

If $\widetilde{G}$ is a $U(\mathbb{Q})$-gain cycle, then $\widetilde{G}=\widetilde{C}_{g}$ is of Type $1$.

Assume that $\widetilde{G}$ is not a $U(\mathbb{Q})$-gain cycle. Combining $r(\widetilde{G})=r(\widetilde{C}_{g})=g-2$ and Lemma \ref{le:2.8}, $V(\widetilde{C}_{g})$ is a dominating set of $\widetilde{G}$. If there exists $v\in V(\widetilde{G})\backslash V(\widetilde{C}_{g})$ such that $|N_{\widetilde{C}_{g}}(v)|=1$, by Lemmas \ref{le:2.6} and \ref{le:2.1}, then
$$r(\widetilde{C}_{g}+v)=r(\widetilde{P}_{g-1})+2=(g-2)+2=g>r(\widetilde{G}),$$
a contradiction.
So $|N_{\widetilde{C}_{g}}(v)|\geq2$ for every $v\in V(\widetilde{G})\backslash V(\widetilde{C}_{g})$.
Combining $g$ is even and Lemma \ref{le:2.9}, $g=4$.
Hence $r(\widetilde{G})=2$ and $|N_{\widetilde{C}_{g}}(v)|=2$ for every $v\in V(\widetilde{G})\backslash V(\widetilde{C}_{g})$.

Let $V(\widetilde{C}_{g})=\{v_{1},v_{2},v_{3},v_{4}\}$ and let $ V(\widetilde{G})\backslash V(\widetilde{C}_{g})=\{u_{1},u_{2},\ldots,u_{n-4}\}$.
Since $g=4$, $N_{\widetilde{C}_{g}}(u_{i})=\{v_{1},v_{3}\}$ or $\{v_{2},v_{4}\}$ ($i=1,2,\ldots,n-4$).
If there exists two different vertices $u_{i},u_{j}$ ($i,j\in\{1,2,\ldots,n-4\}$) such that $N_{\widetilde{C}_{g}}(u_{i})\neq N_{\widetilde{C}_{g}}(u_{j})$, then without loss of generality assumption that $N_{\widetilde{C}_{g}}(u_{i})=\{v_{2},v_{4}\}$ and $N_{\widetilde{C}_{g}}(u_{j})=\{v_{1},v_{3}\}$.
Let $\widetilde{C}_{1}$ be the $U(\mathbb{Q})$-gain cycle with vertices $v_{1},v_{4},v_{3},u_{j}$ in turn.
If $N_{\widetilde{C}_{1}}(u_{i})=\{v_{4}\}$, by Lemma \ref{le:2.6}, then
$r(\widetilde{C}_{1}+u_{i})=4>r(\widetilde{G})$, a contradiction.
Combined with $g=4$, $N_{\widetilde{C}_{1}}(u_{i})=\{v_{4}, u_{j}\}$.
Thus $\widetilde{G}$ is a $U(\mathbb{Q})$-gain complete bipartite graph $\widetilde{K}_{a,b}$, where $a=|N_{\widetilde{G}}(v_{2})|\geq2$ and $b=|N_{\widetilde{G}}(v_{1})|\geq2$.
By Lemma \ref{le:3.1}, all the $\widetilde{C}_{4}$ in $\widetilde{K}_{a,b}$ are of Type $1$.
\quad $\square$\\

\begin{figure}[htbp]
  \centering
  \includegraphics[scale=0.8]{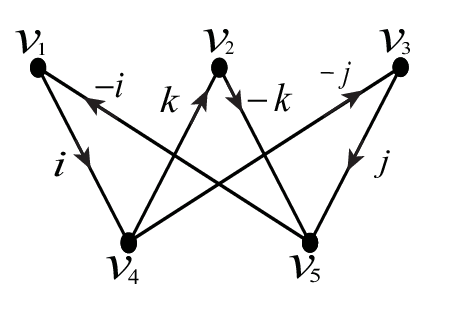}
  \caption{$\widetilde{G}$.}
\end{figure}

\noindent\begin{example}
Consider the $U(\mathbb{Q})$-gain graph $\widetilde{G}=\widetilde{K}_{3,2}$ in Fig. $2$, then
$$\varphi(v_{1}v_{4}v_{2}v_{5}v_{1})=ik(-k)(-i)=(-j)(-k)(-i)=i(-i)=1;$$
$$\varphi(v_{1}v_{4}v_{3}v_{5}v_{1})=i(-j)j(-i)=(-k)j(-i)=i(-i)=1;$$
$$\varphi(v_{2}v_{4}v_{3}v_{5}v_{2})=(-k)(-j)jk=(-i)jk=(-k)k=1.$$
By Definition \ref{de:2.2}, all the $\widetilde{C}_{4}$ in $\widetilde{K}_{3,2}$ are of Type $1$.
Combining $g=4$ and Lemma \ref{le:3.1}, $r(\widetilde{G})=2=g-2$.
Thus $\widetilde{G}$ is an extremal graph which satisfies the condition in Theorem \ref{th:3.2}.
\end{example}

\section{$U(\mathbb{Q})$-gain graphs $\widetilde{G}$ with $r(\widetilde{G})=g-1$}

In this section, we will characterize the $U(\mathbb{Q})$-gain graphs with $r(\widetilde{G})=g-1$. First, we will characterize the rank of $U(\mathbb{Q})$-gain connected bicyclic graphs.

\noindent\begin{lemma}\label{le:4.1}
Let $\widetilde{G}$ be a $U(\mathbb{Q})$-gain $\widetilde{\infty}(p,l,q)$-graph of order $n$ with pendant vertices. Then
\begin{align*}
r(\widetilde{G})\geq\left\{\begin{array}{ll}
p+q,  ~~~~~~~if~p,~q~are~odd;\\
p+q-2, ~~if~p,~q~are~even;\\
p+q-1, ~~otherwise.\\
\end{array}\right.
\end{align*}
This bound is sharp.
\end{lemma}
\noindent\textbf{Proof.}
By Lemmas \ref{le:2.10.}, \ref{le:2.11} and \ref{le:2.12}, we only need to consider the following three cases.

\textbf{Case 1.} $\widetilde{G}$ is obtained from $\widetilde{\infty}(p,1,q)$ by adding $n-p-q+1$ pendant vertices to $x$ which is the common vertex of $\widetilde{C}_{p}$ and $\widetilde{C}_{q}$.

By Lemma \ref{le:2.6}, $r(\widetilde{G})=r(\widetilde{P}_{p-1})+r(\widetilde{P}_{q-1})+2.$
If $p,q$ are odd, by Lemma \ref{le:2.1}, then $r(\widetilde{G})=p-1+q-1+2=p+q$.
If $p,q$ are even, by Lemma \ref{le:2.1}, then $r(\widetilde{G})=p-2+q-2+2=p+q-2$.
Otherwise, by Lemma \ref{le:2.1}, $r(\widetilde{G})=p+q-1-2+2=p+q-1$.
So $r(\widetilde{G})$ attains the above lower bounds.

\textbf{Case 2.} $\widetilde{G}$ is obtained from $\widetilde{\infty}(p,1,q)$ by adding $n-p-q+1$ pendant vertices to a vertex (different from $x$) of $\widetilde{C}_{p}$ or $\widetilde{C}_{q}$.

Without loss of generality, we add pendant vertices to a vertex of $\widetilde{C}_{p}$.
If $p,q$ are odd, by Lemmas \ref{le:2.1} and \ref{le:2.6}, then $r(\widetilde{G})=2+p-1+r(\widetilde{P}_{q-1})= 2+p-1+q-1=p+q$.
If $p,q$ are even, by Lemmas \ref{le:2.1}, \ref{le:2.4} and \ref{le:2.6}, then $r(\widetilde{G})\geq2+p-2+r(\widetilde{P}_{q-1})=2+p-2+q-2=p+q-2$.
If $p$ is odd and $q$ is even, by Lemmas \ref{le:2.1} and \ref{le:2.6}, then $r(\widetilde{G})=2+p-1+r(\widetilde{P}_{q-1})=2+p-1+q-2=p+q-1$.
If $p$ is even and $q$ is odd, by Lemmas \ref{le:2.1}, \ref{le:2.4} and \ref{le:2.6}, then $r(\widetilde{G})\geq2+p-2+r(\widetilde{P}_{q-1})=2+p-2+q-1=p+q-1$.

\textbf{Case 3.} $\widetilde{G}$ is obtained from $\widetilde{\infty}(p,2,q)$ by adding $n-p-q$ pendant vertices to a vertex $y$ of $\widetilde{C}_{p}$ or $\widetilde{C}_{q}$. Without loss of generality, $y\in V(\widetilde{C}_{p})$.

If $p,q$ are odd, by Lemmas \ref{le:2.3} and \ref{le:2.6}, then $r(\widetilde{G})=2+p-1+r(\widetilde{C}_{q})\geq2+p-1+q-1=p+q$.
If $p,q$ are even, by Lemmas \ref{le:2.3}, \ref{le:2.4} and \ref{le:2.6}, then $r(\widetilde{G})\geq2+p-2+r(\widetilde{C}_{q})\geq2+p-2+q-2=p+q-2$.
If $p$ is odd and $q$ is even, by Lemmas \ref{le:2.3} and \ref{le:2.6}, then $r(\widetilde{G})=2+p-1+r(\widetilde{C}_{q})\geq 2+p-1+q-2=p+q-1$.
If $p$ is even and $q$ is odd, by Lemmas \ref{le:2.3}, \ref{le:2.4} and \ref{le:2.6}, then $r(\widetilde{G})\geq2+p-2+r(\widetilde{C}_{q})\geq2+p-2+q-1=p+q-1$.
\quad $\square$

\noindent\begin{lemma}\label{le:4.2}
Let $\widetilde{G}$ be a $U(\mathbb{Q})$-gain $\widetilde{\theta}(p,l,q)$-graph of order $n$ with pendant vertices. If $plq\neq0$, then
\begin{align*}
r(\widetilde{G})\geq\left\{\begin{array}{ll}
p+q+l+1,  ~~if~p~is~odd;\\
p+q+l+2, ~~if~p~is~even.\\
\end{array}\right.
\end{align*}
This bound is sharp.
\end{lemma}
\noindent\textbf{Proof.}
By Lemmas \ref{le:2.10.}, \ref{le:2.11} and \ref{le:2.12}, we only need to consider the following two cases.

\textbf{Case 1.} Adding all pendant vertices to a $3$-degree vertex $x$ of $\widetilde{\theta}(p,l,q)$.

By Lemmas \ref{le:2.6} and \ref{le:2.1},
\begin{align*}
r(\widetilde{G})=\left\{\begin{array}{ll}
p+3+r(\widetilde{P}_{l})+r(\widetilde{P}_{q})\geq p+3+l-1+q-1=p+l+q+1,~if~p~is~odd;\\
p+2+r(\widetilde{P}_{l+q+1})\geq p+2+l+q=p+l+q+2,~~~~~~~~~~~~~~~~if~p~is~even.\\
\end{array}\right.
\end{align*}
If $p$, $l$ and $q$ are odd, or $p$ and $l+q$ are even, then the equality holds.
Thus $r(\widetilde{G})$ attains the above lower bounds.

\textbf{Case 2.} Adding all pendant vertices to a $2$-degree vertex $y$ of $\widetilde{\theta}(p,l,q)$.

Without loss of generality, we assume $y\in V(\widetilde{P}_{p})$.
If $p$ is odd, by Lemmas \ref{le:2.6}, \ref{le:2.1} and \ref{le:2.3}, then
$r(\widetilde{G})=p+3+r(\widetilde{P}_{l})+r(\widetilde{P}_{q})\geq p+3+l-1+q-1=p+l+q+1$,
or $r(\widetilde{G})=p+1+r(\widetilde{C}_{l+q+2})\geq p+1+l+q=p+l+q+1$.
If $p$ is even, by Lemmas \ref{le:2.6} and \ref{le:2.1}, then
$r(\widetilde{G})=p+2+r(\widetilde{P}_{l+q+1})\geq p+2+l+q=p+l+q+2$.
\quad $\square$\\

Next, we consider the case that one of $p,l,q$ is equal to $0$. Without loss of generality, we assume $p=0$.
Similar to the proof for  Lemma \ref{le:4.2}, we can prove the following lemma.

\noindent\begin{lemma}\label{le:4.3}
Let $\widetilde{G}$ be a $U(\mathbb{Q})$-gain $\widetilde{\theta}(0,l,q)$-graph of order $n$ with pendant vertices. Then
\begin{align*}
r(\widetilde{G})\geq\left\{\begin{array}{ll}
l+q+1,  ~~if~l+q~is~odd;\\
l+q+2, ~~if~l+q~is~even.\\
\end{array}\right.
\end{align*}
This bound is sharp.
\end{lemma}

\begin{figure}[htbp]
  \centering
  \includegraphics[scale=0.85]{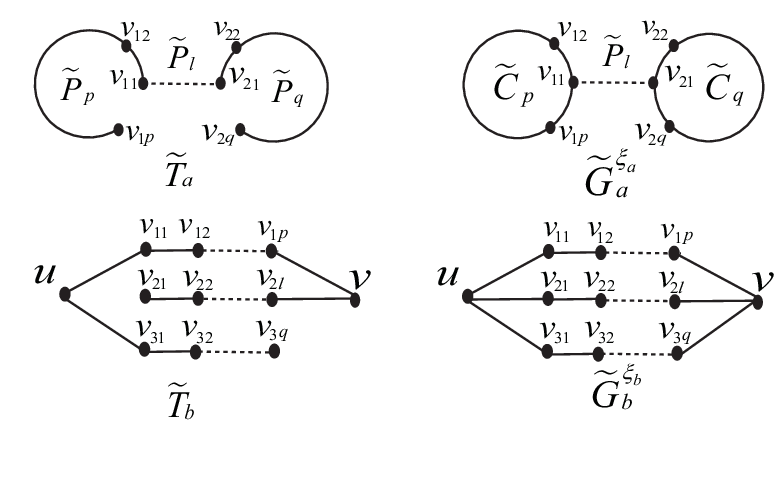}
\caption{ $\widetilde{T}_{a}$, $\widetilde{T}_{b}$, $\widetilde{G}_{a}^{\xi_{a}}$ and $\widetilde{G}_{b}^{\xi_{b}}$.}
\end{figure}

Now, we will characterize the $U(\mathbb{Q})$-gain connected bicyclic graphs with rank $2$, $3$ and $4$.

A \emph{spanning tree} of a connected graph is a subgraph that is a tree with the same vertices as the graph.
Assume $\widetilde{G}$ is connected. Let $\widetilde{T}$ be a spanning tree of $\widetilde{G}$ and let $u$ be a root vertex.
For $v\in V(\widetilde{G})\backslash\{u\}$, let $\widetilde{P}_{vu}$ be the unique path from $v$ to $u$ in $\widetilde{T}$.
Let $\xi(v)=\varphi(\widetilde{P}_{vu})$ and $\xi(u)=1$.
If $xy\in E(\widetilde{G})$ and $d_{\widetilde{G}}(x,u)>d_{\widetilde{G}}(y,u)$, then
$$\varphi^{\xi}(xy)=\xi(x)^{-1}\varphi(xy)\xi(y)
=\varphi(\widetilde{P}_{xu})^{-1}\varphi(xy)\varphi(\widetilde{P}_{yu})
=\varphi(\widetilde{P}_{xu})^{-1}\varphi(\widetilde{P}_{xu})=1.$$
Thus the switching function $\xi$ reduces all the gains on $\widetilde{T}$ to $1$.

Let $\widetilde{T}_{a}$ (resp., $\widetilde{T}_{b}$) be a spanning tree of $\widetilde{G}_{a}$ (resp., $\widetilde{G}_{b}$).
We choose $v_{11}$ (resp., $u$) to be the root vertex of $\widetilde{T}_{a}$ (resp., $\widetilde{T}_{b}$).
After using the above switching function, we have that all edges in $\widetilde{T}_{a}$ (resp., $\widetilde{T}_{b}$) have gain $1$.
So in $\widetilde{G}_{a}^{\xi_{a}}$ (resp., $\widetilde{G}_{b}^{\xi_{b}}$), the gains of the two edges $v_{11}v_{1p}$, $v_{21}v_{2q}$ (resp., $uv_{21}$, $vv_{3q}$) equal to the gains of the two cycles in $\widetilde{G}_{a}$ (resp., $\widetilde{G}_{b}$) and the gains of the other edges have gain $1$.

\begin{figure}[htbp]
  \centering
  \includegraphics[scale=0.6]{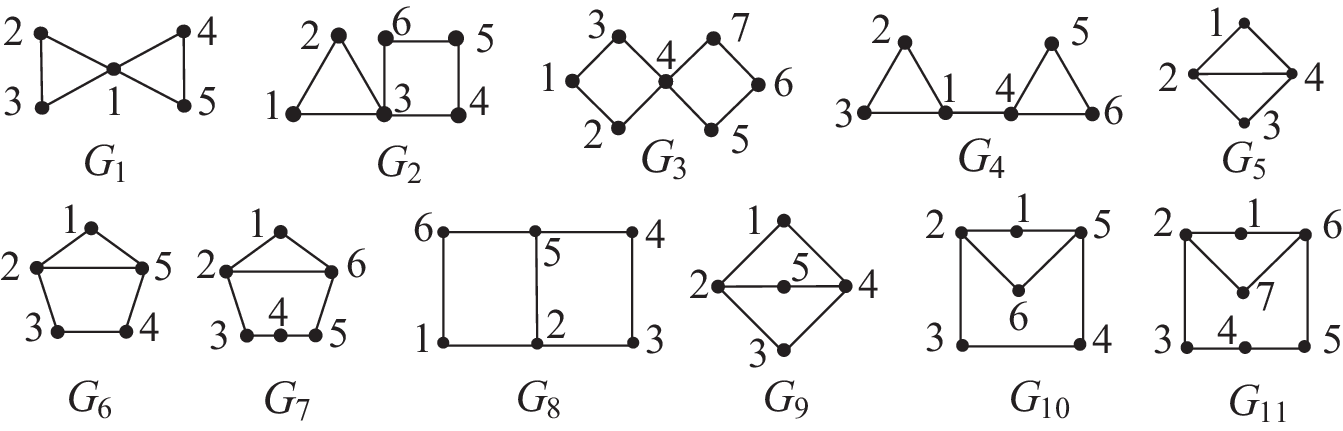}
\caption{The eleven graphs in Theorem 4.4.}
\end{figure}

\noindent\begin{theorem}\label{th:4.4}
Let $\widetilde{G}$ be a $U(\mathbb{Q})$-gain connected bicyclic graph without pendant vertices. Then the
cases of graphs with rank $2$, $3$ and $4$ are given in Table $1$.
\end{theorem}

 \begin{minipage}[hbt]{1\columnwidth}
\textbf{Table 1}

 {\small { The gain conditions for each gain graph in Theorem \ref{th:4.4} with rank 2, 3 or 4.}}
  \vskip1mm
 {
\begin{tabular}{p{0.5cm}p{12cm}}
\hline
 $r(\widetilde{G})$& Gain graph $\widetilde{G}$ and its gain conditions \\
\hline
$2$

&$\widetilde{G}_{5}$.~~The subgraph induced on vertices $1,~2,~4$ is of Type $4$ and the subgraph induced on vertices $1,~2,~3,~4$ is of Type $1$.\\

&$\widetilde{G}_{9}$.~~The subgraphs induced on vertices $1,~2,~5,~4$ and vertices $1,~2,~3,~4$ are of Type $1$.
\\
\hline
$3$ & $\widetilde{G}_{5}$.~~The subgraph induced on vertices $1,~2,~4$ is of Type $3$ and the subgraph induced on vertices $1,~2,~3,~4$ is of Type $1$.
\\
\hline
$4$ & $\widetilde{G}_{1}$.~ $Re\left(\varphi(v_{1}v_{3}v_{2}v_{1})\right)+Re\left(\varphi(v_{1}v_{5}v_{4}v_{1})\right)=0$.
\\
&$\widetilde{G}_{2}$.~~The subgraph induced on vertices $3,~4,~5,~6$ is of Type $1$ and the subgraph induced on vertices $1,~2,~3$ is of Type $4$.\\

&$\widetilde{G}_{3}$.~~The subgraphs induced on vertices $1,~2,~4,~3$ and vertices $4,~5,~6,~7$ are of Type $1$.\\

&$\widetilde{G}_{5}$.~~The subgraph induced on vertices $1,~2,~4$ is of Type $3$ or $4$ and the subgraph induced on vertices $1,~2,~3,~4$ is of Type $2$.\\

&$\widetilde{G}_{6}$.~
$Re\left(\varphi(v_{1}v_{2}v_{5}v_{1})\right)-Re\left(\varphi(v_{1}v_{2}v_{3}v_{4}v_{5}v_{1})\right)=0$.\\

&$\widetilde{G}_{7}$.~~The subgraph induced on vertices $1,~2,~6$ is of Type $4$ and the subgraph induced on vertices $1,~2,~3,~4,~5,~6$ is of Type $1$.\\

&$\widetilde{G}_{8}$.~
$\varphi(v_{1}v_{2}v_{3}v_{4}v_{5}v_{6}v_{1})-\varphi(v_{1}v_{2}v_{5}v_{6}v_{1})+1=0$.\\

&$\widetilde{G}_{9}$.~~The subgraph induced on vertices $1,~2,~3,~4$ is of Type $1$ and the subgraph induced on vertices $1,~2,~5,~4$ is of Type $2$, or the subgraph induced on vertices $1,~2,~3,~4$ is of Type $2$ and the subgraph induced on vertices $1,~2,~5,~4$ is of Type $1$ or $2$.\\

&$\widetilde{G}_{10}$.~~The subgraph induced on vertices $1,~2,~6,~5$ is of Type $1$ and the subgraph induced on vertices $1,~2,~3,~4,~5$ is of Type $4$.\\

&$\widetilde{G}_{11}$.~~The subgraphs induced on vertices $1,~2,~7,~6$ and vertices $1,~2,~3,~4,~5,~6$ are of Type $1$.
\\
\hline
\end{tabular} }
\end{minipage}

\noindent\textbf{Proof.}
Let $G$ be the underlying graph of $\widetilde{G}$.

\textbf{Case 1.} $\widetilde{G}=\widetilde{\infty}(p,l,q)$. The switching function $\xi_{1}$ and the spanning tree in $\widetilde{G}^{\xi_{1}}$ is similar to those in the above graph $\widetilde{G}_{a}^{\xi_{a}}$ in Fig. 3.

Without loss of generality, $p\leq q$.
If $(p,l,q)\in A=\{(3,1,3),(3,1,4),(4,1,4),(3,2,3)\}$, then $G$ is one of the graphs $G_{1}-G_{4}$ in Fig. 4.

If $G=G_{1}$, then $\varphi^{\xi_{1}}(e_{13})=\varphi(v_{1}v_{3}v_{2}v_{1})$ and $\varphi^{\xi_{1}}(e_{15})=\varphi(v_{1}v_{5}v_{4}v_{1})$.
Let $a=\varphi^{\xi_{1}}(e_{13})$ and $b=\varphi^{\xi_{1}}(e_{15})$.
After calculation, we get $r(\widetilde{G}^{\xi_{1}})=4+r(a+\overline{a}+b+\overline{b})$.
So $r(\widetilde{G}^{\xi_{1}})=4$ if and only if $a+\overline{a}+b+\overline{b}=0$, i.e., $Re(\varphi(v_{1}v_{3}v_{2}v_{1}))+Re(\varphi(v_{1}v_{5}v_{4}v_{1}))=0$.

If $G=G_{2}$, by Lemmas \ref{le:2.3} and \ref{le:2.13}, then $r(\widetilde{G})=4$ if and only if the subgraph induced on vertices $3,4,5,6$ is of Type $1$ and the subgraph induced on vertices $1,2,3$ is of Type $4$.

If $G=G_{3}$, by Lemmas \ref{le:2.3} and \ref{le:2.13}, then $r(\widetilde{G})=4$ if and only if the subgraphs induced on vertices $1,2,4,3$ and vertices $4,5,6,7$ are of Type $1$.

If $G=G_{4}$, then $\varphi^{\xi_{1}}(e_{13})=\varphi(v_{1}v_{3}v_{2}v_{1})$ and $\varphi^{\xi_{1}}(e_{46})=\varphi(v_{4}v_{6}v_{5}v_{4})$.
Let $a=\varphi^{\xi_{1}}(e_{13})$ and $b=\varphi^{\xi_{1}}(e_{46})$.
After calculation, we have
\begin{align*}
 r(\widetilde{G}^{\xi_{1}})=4+r\left
 (\begin{array}{cccccccc}
 -a-\overline{a} & 1\\
 -1 & b+\overline{b}\\
\end{array}
 \right),
 \end{align*}
so $r(\widetilde{G}^{\xi_{1}})>4$.

If $(p,l,q)\notin A$, by Lemmas \ref{le:2.1}, \ref{le:2.3} and \ref{le:2.13}, then $r(\widetilde{G})>4$.

\textbf{Case 2.} $\widetilde{G}=\widetilde{\theta}(p,l,q)$. The switching function $\xi_{2}$ and the spanning tree in $\widetilde{G}^{\xi_{2}}$ is similar to those in the above graph $\widetilde{G}_{b}^{\xi_{b}}$ in Fig. 3.

Without loss of generality, $p\leq l\leq q$.
If $(p,l,q)\in B=\{(0,1,1),(0,1,2),(0,1,3),(0,2,2),$ $(1,1,1)$, $(1,1,2),(1,1,3)\}$, then $G$ is one of the graphs $G_{5}-G_{11}$ in Fig. 4.

If $G=G_{5}$, then $\varphi^{\xi_{2}}(e_{24})=\varphi(v_{1}v_{2}v_{4}v_{1})$ and $\varphi^{\xi_{1}}(e_{34})=\varphi(v_{1}v_{2}v_{3}v_{4}v_{1})$.
Let $a=\varphi^{\xi_{2}}(e_{24})$ and $b=\varphi^{\xi_{2}}(e_{34})$.
After calculation, we get
\begin{align*}
 r(\widetilde{G}^{\xi_{2}})=2+r\left
 (\begin{array}{cccccccc}
 0 & b-1\\
 \overline{b}-1 & -a-\overline{a}\\
\end{array}
 \right),
 \end{align*}
so
\begin{align*}
r(\widetilde{G}^{\xi_{2}})=\left\{\begin{array}{ll}
2, & if~a+\overline{a}=0,~b=1;\\
3,& if~a+\overline{a}\neq0,~b=1;\\
4,& if~b\neq1,\\
\end{array}\right.
\end{align*}
i.e., $r(\widetilde{G}^{\xi_{2}})=2$ if and only if the subgraph induced on vertices $1,2,4$ is of Type $4$ and the subgraph induced on vertices $1,2,3,4$ is of Type $1$;
$r(\widetilde{G}^{\xi_{2}})=3$ if and only if the subgraph induced on vertices $1,2,4$ is of Type $3$ and the subgraph induced on vertices $1,2,3,4$ is of Type $1$;
$r(\widetilde{G}^{\xi_{2}})=4$ if and only if the subgraph induced on vertices $1,2,4$ is of Type $3$ or $4$ and the subgraph induced on vertices $1,2,3,4$ is of Type $2$.

If $G=G_{6}$, then $\varphi^{\xi_{2}}(e_{25})=\varphi(v_{1}v_{2}v_{5}v_{1})$ and $\varphi^{\xi_{2}}(e_{45})=\varphi(v_{1}v_{2}v_{3}v_{4}v_{5}v_{1})$.
Let $a=\varphi^{\xi_{2}}(e_{25})$ and $b=\varphi^{\xi_{2}}(e_{45})$.
After calculation, we get $r(\widetilde{G}^{\xi_{2}})=4+r(-a-\overline{a}+b+\overline{b})$.
So $r(\widetilde{G}^{\xi_{2}})=4$ if and only if $-a-\overline{a}+b+\overline{b}=0$, i.e., $Re(\varphi(v_{1}v_{2}v_{5}v_{1}))-Re(\varphi(v_{1}v_{2}v_{3}v_{4}v_{5}v_{1}))=0$.

If $G=G_{7}$, then $\varphi^{\xi_{2}}(e_{26})=\varphi(v_{1}v_{2}v_{6}v_{1})$ and $\varphi^{\xi_{2}}(e_{56})=\varphi(v_{1}v_{2}v_{3}v_{4}v_{5}v_{6}v_{1})$.
Let $a=\varphi^{\xi_{2}}(e_{26})$ and $b=\varphi^{\xi_{2}}(e_{56})$.
After calculation, we get
\begin{align*}
 r(\widetilde{G}^{\xi_{2}})=4+r\left
 (\begin{array}{cccccccc}
 0 & b+1\\
 \overline{b}+1 & -a-\overline{a}\\
\end{array}
 \right),
 \end{align*}
so $r(\widetilde{G}^{\xi_{2}})=4$ if and only if $a+\overline{a}=0$ and $b=-1$, i.e., the subgraph induced on vertices $1,2,6$ is of Type $4$ and the subgraph induced on vertices $1,2,3,4,5,6$ is of Type $1$.

If $G=G_{8}$, then $\varphi^{\xi_{2}}(e_{25})=\varphi(v_{1}v_{2}v_{5}v_{6}v_{1})$ and $\varphi^{\xi_{2}}(e_{45})=\varphi(v_{1}v_{2}v_{3}v_{4}v_{5}v_{6}v_{1})$.
Let $a=\varphi^{\xi_{2}}(e_{25})$ and $b=\varphi^{\xi_{2}}(e_{45})$.
After calculation, we get
\begin{align*}
 r(\widetilde{G}^{\xi_{2}})=4+r\left
 (\begin{array}{cccccccc}
 0 & \overline{b}-\overline{a}+1\\
 b-a+1 & 0\\
\end{array}
 \right),
 \end{align*}
so $r(\widetilde{G}^{\xi_{2}})=4$ if and only if $b-a+1=0$, i.e., $\varphi(v_{1}v_{2}v_{3}v_{4}v_{5}v_{6}v_{1})-\varphi(v_{1}v_{2}v_{5}v_{6}v_{1})+1=0$.

If $G=G_{9}$, then $\varphi^{\xi_{2}}(e_{25})=\varphi(v_{1}v_{2}v_{5}v_{4}v_{1})$ and $\varphi^{\xi_{2}}(e_{34})=\varphi(v_{1}v_{2}v_{3}v_{4}v_{1})$.
Let $a=\varphi^{\xi_{2}}(e_{25})$ and $b=\varphi^{\xi_{2}}(e_{34})$.
After calculation, we get
\begin{align*}
 r(\widetilde{G}^{\xi_{2}})=2+r\left
 (\begin{array}{cccccccc}
 0 & b-1 & 0\\
 \overline{b}-1 & 0 & 1-a\\
 0 & 1-\overline{a} & 0 \\
\end{array}
 \right),
 \end{align*}
so
\begin{align*}
r(\widetilde{G}^{\xi_{2}})=\left\{\begin{array}{ll}
2, & if~a=1,~b=1;\\
4,& otherwise,\\
\end{array}\right.
\end{align*}
i.e., $r(\widetilde{G}^{\xi_{2}})=2$ if and only if the subgraphs induced on vertices $1,2,3,4$ and vertices $1,2,5,4$ are of Type $1$;
$r(\widetilde{G}^{\xi_{2}})=4$ if and only if the subgraph induced on vertices $1,2,3,4$ is of Type $1$ and the subgraph induced on vertices $1,2,5,4$ is of Type $2$, or the subgraph induced on vertices $1,2,3,4$ is of Type $2$ and the subgraph induced on vertices $1,2,5,4$ is of Type $1$ or $2$.

If $G=G_{10}$, then $\varphi^{\xi_{2}}(e_{26})=\varphi(v_{1}v_{2}v_{6}v_{5}v_{1})$ and $\varphi^{\xi_{2}}(e_{45})=\varphi(v_{1}v_{2}v_{3}v_{4}v_{5}v_{1})$.
Let $a=\varphi^{\xi_{2}}(e_{26})$ and $b=\varphi^{\xi_{2}}(e_{45})$.
After calculation,
\begin{align*}
 r(\widetilde{G}^{\xi_{2}})=4+r\left
 (\begin{array}{cccccccc}
 b+\overline{b} & 1-a\\
 1-\overline{a} & 0\\
\end{array}
 \right),
 \end{align*}
so $r(\widetilde{G}^{\xi_{2}})=4$ if and only if $a=1$ and $b+\overline{b}=0$, i.e., the subgraph induced on vertices $1,2,6,5$ is of Type $1$ and the subgraph induced on vertices $1,2,3,4,5$ is of Type $4$.

If $G=G_{11}$, then $\varphi^{\xi_{2}}(e_{27})=\varphi(v_{1}v_{2}v_{7}v_{6}v_{1})$ and $\varphi^{\xi_{2}}(e_{56})=\varphi(v_{1}v_{2}v_{3}v_{4}v_{5}v_{6}v_{1})$.
Let $a=\varphi^{\xi_{2}}(e_{27})$ and $b=\varphi^{\xi_{2}}(e_{56})$.
After calculation,
\begin{align*}
 r(\widetilde{G}^{\xi_{2}})=4+r\left
 (\begin{array}{cccccccc}
 0 & b+1 & 0\\
 \overline{b}+1 & 0 & 1-a\\
 0 & 1-\overline{a} & 0\\
\end{array}
 \right),
 \end{align*}
so $r(\widetilde{G}^{\xi_{2}})=4$ if and only if $a=1$ and $b=-1$, i.e., the subgraphs induced on vertices $1,2,7,6$ and vertices $1,2,3,4,5,6$ are of Type $1$.

If $(p,l,q)\notin B$, by Lemmas \ref{le:2.1}, \ref{le:2.3} and \ref{le:2.13}, then $r(\widetilde{G})>4$.
\quad $\square$

\begin{figure}[htbp]
\begin{minipage}[hbt]{1\columnwidth}
\textbf{Table 2}

 {\small {The gain conditions for each gain graph in Theorem \ref{th:4.5} satisfying $r(\widetilde{G})=4$.}}
  \vskip1mm
 {
\begin{tabular}{p{5cm}p{8.8cm}}
\hline
Gain graphs $\widetilde{G}$ & Gain conditions of $\widetilde{G}$\\
\hline
$\widetilde{G}_{12}$ & The subgraph induced on vertices $1,~2,~3$ is of Type $4$. \\

\hline
$\widetilde{G}_{13},~\widetilde{G}_{14},~\widetilde{G}_{19},~\widetilde{G}_{20}$ & Any gain. \\
\hline
$\widetilde{G}_{15},~\widetilde{G}_{16}$ & The subgraph induced on vertices $1,~2,~4$ is of Type $4$ and the subgraph induced on vertices $1,~2,~3,~4$ is of Type $1$. \\
\hline
$\widetilde{G}_{17}$ & The subgraph induced on vertices $1,~2,~3,~4$ is of Type $1$. \\
\hline
$\widetilde{G}_{18}$ & The subgraph induced on vertices $1,~2,~3,~4$ is of Type $1$. \\
\hline
$\widetilde{G}_{21},~\widetilde{G}_{22}$ & The subgraphs induced on vertices $1,~2,~5,~4$ and vertices $1,~2,~3,~4$ are of Type $1$. \\

\hline
\end{tabular} }
\end{minipage}
\end{figure}

\begin{figure}[htbp]
  \centering
  \includegraphics[scale=0.6]{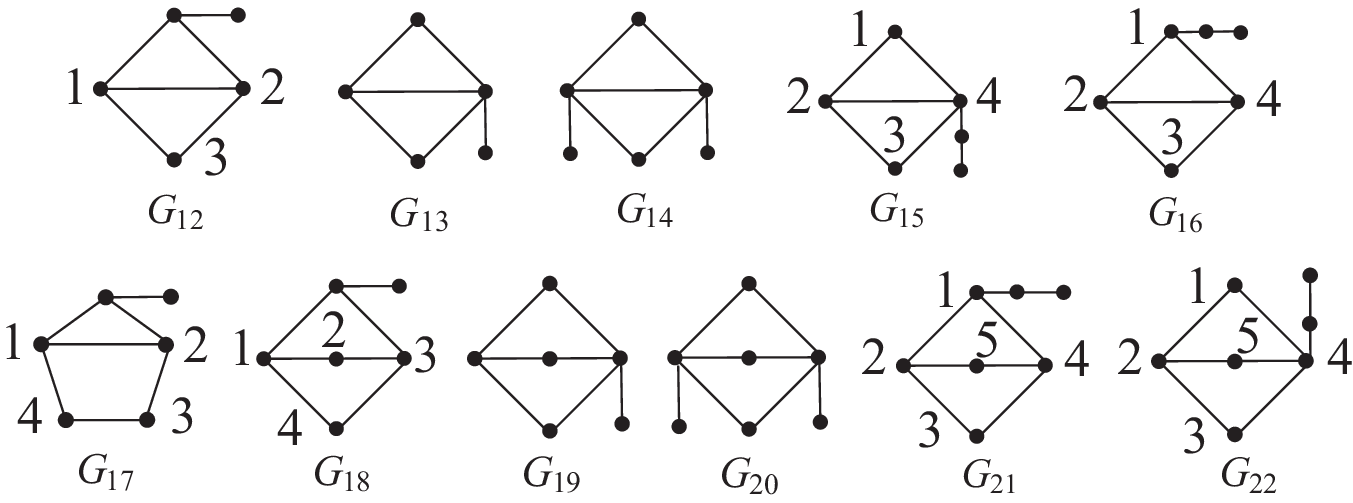}
\caption{The eleven graphs in Theorem 4.5.}
\end{figure}

\noindent\begin{theorem}\label{th:4.5}
Let $\widetilde{G}$ be a $U(\mathbb{Q})$-gain connected bicyclic graph with pendant vertices but no pendant twins. Then the cases of graphs with rank $4$ are given in Table $2$.
\end{theorem}
\noindent\textbf{Proof.}
By Lemmas \ref{le:4.1}, \ref{le:4.2}, \ref{le:4.3} and $r(\widetilde{G})=4$, we have that $\widetilde{G}$ is a $\widetilde{\theta}(p,l,q)$-graph and $(p,l,q)\in\{(0,1,1),(0,1,2),(1,$ $1,1)\}$. Let $G$ be the underlying graph of $\widetilde{G}$.

\textbf{Case 1.} $(p,l,q)=(0,1,1)$.

\textbf{Subcase 1.1.} $\widetilde{G}-\widetilde{\theta}(0,1,1)$ is some isolated vertices.

If $|V(\widetilde{G}-\widetilde{\theta}(0,1,1))|=1$, then $G$ is $G_{12}$ or $G_{13}$ (see Fig. 5).
When $G=G_{12}$, by Lemmas \ref{le:2.3} and \ref{le:2.6}, $r(\widetilde{G})=4$ if and only if the subgraph induced on vertices $1,2,3$ is of Type $4$.
When $G=G_{13}$, by Lemma \ref{le:2.6}, $r(\widetilde{G})=4$.

If $|V(\widetilde{G}-\widetilde{\theta}(0,1,1))|=2$, by Lemma \ref{le:2.6}, then $r(\widetilde{G})=4$ if and only if $G=G_{14}$ (see Fig. 5), otherwise $r(\widetilde{G})=6$.

If $|V(\widetilde{G}-\widetilde{\theta}(0,1,1))|\geq 3$, by Lemma \ref{le:2.6}, then $r(\widetilde{G})\geq 6$.

\textbf{Subcase 1.2.} $\widetilde{G}-\widetilde{\theta}(0,1,1)=\widetilde{P}_{2}$.

In this subcase, $G$ is $G_{15}$ or $G_{16}$ (see Fig. 5). By Lemma \ref{le:2.6} and Theorem \ref{th:4.4}, $r(\widetilde{G})=4$ if and only if the subgraph induced on vertices $1,~2,~4$ is of Type $4$ and the subgraph induced on vertices $1,~2,~3,~4$ is of Type $1$.

\textbf{Subcase 1.3.} $\widetilde{G}-\widetilde{\theta}(0,1,1)$  contains the union of $\widetilde{P}_{2}$ and isolated vertices, or contains
$\widetilde{P}_{3}$ as an induced subgraph.

By Lemmas \ref{le:2.6} and \ref{le:2.3}, $r(\widetilde{G})\geq6$.

\textbf{Case 2.} $(p,l,q)=(0,1,2)$.

\textbf{Subcase 2.1.} $\widetilde{G}-\widetilde{\theta}(0,1,2)$ is some isolated vertices.

If $|V(\widetilde{G}-\widetilde{\theta}(0,1,2))|=1$, then $r(\widetilde{G})=4$ if and only if $G=G_{17}$ and the subgraph induced on vertices $1,2,3,4$ is of Type $1$, otherwise $r(\widetilde{G})=6$.

If $|V(\widetilde{G}-\widetilde{\theta}(0,1,2))|\geq2$, by Lemma \ref{le:2.6}, then $r(\widetilde{G})\geq6$.

\textbf{Subcase 2.2.} $\widetilde{G}-\widetilde{\theta}(0,1,2)$ contains $\widetilde{P}_{2}$ as an induced subgraph.

By Lemmas \ref{le:2.6}, \ref{le:2.3} and Theorem \ref{th:4.4}, $r(\widetilde{G})\geq6$.

\textbf{Case 3.} $(p,l,q)=(1,1,1)$.

\textbf{Subcase 3.1.} $\widetilde{G}-\widetilde{\theta}(1,1,1)$ is some isolated vertices.

If $|V(\widetilde{G}-\widetilde{\theta}(1,1,1))|=1$, then $G$ is $G_{18}$ or $G_{19}$ (see Fig. 5).
When $G=G_{18}$, by Lemmas \ref{le:2.3} and \ref{le:2.6}, $r(\widetilde{G})=4$ if and only if the subgraph induced on vertices $1,2,3,4$ is of Type $1$.
When $G=G_{19}$, by Lemma \ref{le:2.6}, $r(\widetilde{G})=4$.

If $|V(\widetilde{G}-\widetilde{\theta}(1,1,1))|=2$, by Lemma \ref{le:2.6}, then $r(\widetilde{G})=4$ if and only if $G=G_{20}$, otherwise $r(\widetilde{G})=6$.

If $|V(\widetilde{G}-\widetilde{\theta}(1,1,1))|\geq 3$, by Lemma \ref{le:2.6}, then $r(\widetilde{G})\geq 6$.

\textbf{Subcase 3.2.} $\widetilde{G}-\widetilde{\theta}(1,1,1)=\widetilde{P}_{2}$.

In this subcase, $G$ is $G_{21}$ or $G_{22}$ (see Fig. 5). By Lemma \ref{le:2.6} and Theorem \ref{th:4.4}, $r(\widetilde{G})=4$ if and only if the subgraphs induced on vertices $1,2,5,4$ and vertices $1,2,3,4$ are of Type $1$.

\textbf{Subcase 3.3.} $\widetilde{G}-\widetilde{\theta}(1,1,1)$  contains the union of $\widetilde{P}_{2}$ and isolated vertices, or contains
$\widetilde{P}_{3}$ as an induced subgraph.

By Lemmas \ref{le:2.6} and \ref{le:4.3}, $r(\widetilde{G})\geq6$.

This proof is complete.\quad $\square$\\

We call distinct vertices $x,y$ of a $U(\mathbb{Q})$-gain graph $\widetilde{G}$ \emph{multiple vertices} of $\widetilde{G}$ if $N_{\widetilde{G}}(x)=N_{\widetilde{G}}(y)$ and there exists $k\in\mathbb{Q}$ such that $k\neq0$ and $\varphi_{xz}=k\varphi_{yz}$ for all $z\in N_{\widetilde{G}}(x)$.
A $U(\mathbb{Q})$-gain graph is said to be a \emph{reduced graph} if it has no multiple vertices.
If $\widetilde{G}$ is a $U(\mathbb{Q})$-gain graph and $u$ is a vertex with a multiple vertex in $\widetilde{G}$,
then removing $u$ from $\widetilde{G}$ is called a reduction of $\widetilde{G}$.
By applying successive reductions to any $U(\mathbb{Q})$-gain graph $\widetilde{G}$, we can obtain a reduced graph of $\widetilde{G}$.

\noindent\begin{lemma}\label{le:4.8}
Let $\widetilde{G}$ be a $U(\mathbb{Q})$-gain graph of order $n$ and let $\widetilde{H}$ be the reduced graph of $\widetilde{G}$. Then $r(\widetilde{G})=r(\widetilde{H})$.
\end{lemma}
\noindent\textbf{Proof.}
Let $x,y$ be multiple vertices of $\widetilde{G}$ and let $\varphi_{xz}=k\varphi_{yz}$ ($k\in\mathbb{Q}$ and $k\neq0$) for all $z\in N_{\widetilde{G}}(x)$.
Let $A(\widetilde{G})=(\alpha_{1},\alpha_{2},\ldots,\alpha_{n})^{T}$.
For convenience, we assume that $\alpha_{1},\alpha_{2}$ are the vectors corresponding to $x,y$ in $A(\widetilde{G})$, respectively.
Then $\alpha_{1}=k\alpha_{2}$. Thus
$$r(\widetilde{G})=r((k\alpha_{2},\alpha_{2},\ldots,\alpha_{n})^{T})=r((\alpha_{2},\ldots,\alpha_{n})^{T})
=r(\widetilde{G}-x).$$
So reductions do not change the rank of $\widetilde{G}$. Thus $r(\widetilde{G})=r(\widetilde{H})$.
\quad $\square$\\

Next, we prove the rank of $U(\mathbb{Q})$-gain complete graphs of order $4$.

\noindent\begin{lemma}\label{le:4.11}
Let $\widetilde{K}_{4}$ be a $U(\mathbb{Q})$-gain complete graph.
Then $r(\widetilde{K}_{4})=4$.
\end{lemma}
\noindent\textbf{Proof.}
Suppose $V(\widetilde{K}_{4})=\{v_{1},v_{2},v_{3},v_{4}\}$.
Let $\xi$ be the switching function such that $\xi(v_{1})=1$ and $\xi(v_{i})=\varphi_{v_{i}v_{1}}$ for $i=2,3,4$.
Thus $\varphi^{\xi}(v_{2}v_{3})=\varphi(v_{1}v_{2}v_{3}v_{1})$, $\varphi^{\xi}(v_{2}v_{4})=\varphi(v_{1}v_{2}v_{4}v_{1})$ and $\varphi^{\xi}(v_{3}v_{4})=\varphi(v_{1}v_{3}v_{4}v_{1})$.
Let $a=\varphi(v_{1}v_{2}v_{3}v_{1})$, $b=\varphi(v_{1}v_{2}v_{4}v_{1})$ and $c=\varphi(v_{1}v_{3}v_{4}v_{1})$.
So
\begin{align*}
 r(\widetilde{K}_{4}^{\xi})=r\left
 (\begin{array}{cccccccc}
 0 & 1 & 1 & 1\\
 1 & 0 & a & b\\
 1 & \overline{a} & 0 & c\\
 1 & \overline{b} & \overline{c} & 0\\
\end{array}
 \right)
 =2+r\left
 (\begin{array}{cccccccc}
 -a-\overline{a} & c-b-\overline{a}\\
 \overline{c}-a-\overline{b} & -b-\overline{b}\\
\end{array}
 \right)=2+r(A_{1}).
 \end{align*}
After calculation, the rows of $A_{1}$ are not left linearly independent.
Thus $r(\widetilde{K}_{4})=4$.
\quad $\square$\\

Now, we will characterize all  $U(\mathbb{Q})$-gain graphs with rank 2.
\noindent\begin{theorem}\label{th:4.10}
Let $\widetilde{G}$ be a connected $U(\mathbb{Q})$-gain graph.
Then $r(\widetilde{G})=2$ if and only if $\widetilde{G}$ is one of the following types:
\begin{enumerate}[(a)]
  \item  $\widetilde{G}=\widetilde{K}_{a,b}$ and all the $\widetilde{C}_{4}$ in $\widetilde{K}_{a,b}$ are of Type 1;
  \item The reduced graph of $\widetilde{G}$ is $\widetilde{C}_{3}$, which is of Type $4$.
\end{enumerate}
\end{theorem}
\noindent\textbf{Proof.}
\textbf{Sufficiency:} If $\widetilde{G}=\widetilde{K}_{a,b}$ and all the $\widetilde{C}_{4}$ in $\widetilde{K}_{a,b}$ are of Type 1, by Lemma \ref{le:3.1}, then $r(\widetilde{G})=2$.
If the reduced graph of $\widetilde{G}$ is $\widetilde{C}_{3}$, which is of Type $4$, then by Lemmas \ref{le:4.8} and \ref{le:2.3}, we have $r(\widetilde{G})=r(\widetilde{C}_{3})=2$.

\textbf{Necessity:}

\textbf{Case 1.} $\widetilde{G}$ is a $U(\mathbb{Q})$-gain bipartite graph.

Let $V(\widetilde{G})=X\cup Y$ and $X\cap Y=\emptyset$.
We claim that $\widetilde{G}$ is a $U(\mathbb{Q})$-gain complete bipartite graph.
On the contrary, there exists two vertices $x_{0}\in X$ and $y_{0}\in Y$ that are not adjacent.
Since $\widetilde{G}$ is connected, we can find a vertex $y'\in Y$ that is adjacent to $x_{0}$
and a vertex $x'\in X$ that is adjacent to $y_{0}$.
By Lemma \ref{le:2.6}, $r(x_{0}+y_{0}+x'+y')=4>r(\widetilde{G})$, a contradiction.
So any $x\in X$ is adjacent to all the vertices $y\in Y$.
Combined with Lemma \ref{le:3.1}, $\widetilde{G}$ is a $U(\mathbb{Q})$-gain complete bipartite graph $\widetilde{K}_{a,b}$ and all the $\widetilde{C}_{4}$ in $\widetilde{K}_{a,b}$ are of Type 1, where $|X|=a$ and $|Y|=b$.

\textbf{Case 2.} $\widetilde{G}$ is not a $U(\mathbb{Q})$-gain bipartite graph.

In this case, $\widetilde{G}$ has an odd cycle $\widetilde{C}_{q}$ as an induced subgraph.
If $q\geq 5$, by Lemma \ref{le:2.3}, then $r(\widetilde{C}_{q})\geq 4>r(\widetilde{G})$, a contradiction.
Therefore, $q=3$.

\textbf{Claim.} $\widetilde{G}$ is a $U(\mathbb{Q})$-gain complete tripartite graph.

Let $x_{0},~y_{0},~z_{0}$ be the tree vertices that induce the $\widetilde{C}_{3}$ of $\widetilde{G}$.
If there exists a vertex $v\in V(\widetilde{G})\backslash V(\widetilde{C}_{3})$, then $|N_{\widetilde{C}_{3}}(v)|=2$.
On the contrary, if $|N_{\widetilde{C}_{3}}(v)|=1$, by Lemma \ref{le:2.6}, then $r(\widetilde{C}_{3}+v)=4>r(\widetilde{G})$, a contradiction.
If $|N_{\widetilde{C}_{3}}(v)|=3$, by Lemma \ref{le:4.11}, then $r(\widetilde{C}_{3}+v)=r(\widetilde{K}_{4})=4>r(\widetilde{G})$, a contradiction.

Let $X$ (resp., $Y$, $Z$) be the set of vertices that are not adjacent to $x_{0}$ (resp., $y_{0}$, $z_{0}$).
If $x'\in X$ is not adjacent to $y'\in Y$, by Lemma \ref{le:2.6}, then $r(y_{0}+z_{0}+x'+y')=4>r(\widetilde{G})$, a contradiction.
Hence, any vertex $x\in X$ is adjacent to all the vertices $y\in Y$.
Similarly, any vertex $x\in X$ is adjacent to all the vertices $z\in Z$ and any vertex $y\in Y$ is adjacent to all the vertices $z\in Z$.
Therefore, $\widetilde{G}$ is a $U(\mathbb{Q})$-gain complete tripartite graph $\widetilde{K}_{r,s,t}$, where $|X|=r,|Y|=s$ and $|Z|=t$.

We may assume that
\begin{align*}
 &A(\widetilde{K}_{r,s,t})=\left (
 \begin{array}{ccccccc}
 \textbf{0} & A_{XY} & A_{XZ}\\
 A_{XY}^{\ast} & \textbf{0} & A_{YZ}\\
 A_{XZ}^{\ast} & A_{YZ}^{\ast} & \textbf{0}\\
 \end{array}
 \right)
 =\left (
 \begin{array}{ccccccc}
 A_{1}\\
 A_{2}\\
 A_{3}\\
 \end{array}
 \right),
  \end{align*}
where $A_{XY}$ is a submatrix representing the connection relationship between vertex set $X$ and vertex set $Y$, and
$A_{XY}^{\ast}$ (resp., $A_{XZ}^{\ast}$, $A_{YZ}^{\ast}$) is the conjugate transpose of $A_{XY}$ (resp., $A_{XZ}$, $A_{YZ}$).
We know that $A_{XY}$, $A_{XZ}$ and $A_{YZ}$ contain nonzero quaternion as entries.
If there exists a row submatrix $A_{i}$ ($i\in\{1,2,3\}$) such that $r(A_{i})\geq 2$, then $r(\widetilde{K}_{r,s,t})\geq 2r(A_{i})\geq 4$, a contradiction.
So, for any $1\leq i\leq3$, each row vector of $A_{i}$ is the nonzero multiple of the first row vector.
It follows that the reduced graph of $\widetilde{K}_{r,s,t}$ is $\widetilde{C}_{3}$.
Combining $r(\widetilde{G})=2$ and Lemma \ref{le:4.8}, $r(\widetilde{C}_{3})=r(\widetilde{G})=2$.
By Lemma \ref{le:2.3}, $\widetilde{C}_{3}$ is of Type $4$.
\quad $\square$\\

Now, we will characterize the $U(\mathbb{Q})$-gain graphs with $r(\widetilde{G})=g-1$.

\noindent\begin{theorem}\label{th:5.1}
Let $\widetilde{G}$ be a connected $U(\mathbb{Q})$-gain graph with girth $g$. Then $r(\widetilde{G})=g-1$ if and only if  $\widetilde{G}$ is one of the following types:
\begin{enumerate}[(a)]
  \item $\widetilde{G}=\widetilde{C}_{g}$ is of Type $4$;
  \item The reduced graph of $\widetilde{G}$ is $\widetilde{C}_{3}$, which is of Type $4$.
\end{enumerate}
\end{theorem}

\noindent\textbf{Proof.}
\textbf{Sufficiency:}
If $\widetilde{G}=\widetilde{C}_{g}$ is of Type $4$, by Lemma \ref{le:2.3}, then $r(\widetilde{G})=g-1$.
If the reduced graph of $\widetilde{G}$ is $\widetilde{C}_{3}$, which is of Type $4$, then by Lemmas \ref{le:4.8} and \ref{le:2.3}, we have $r(\widetilde{G})=r(\widetilde{C}_{3})=2=g-1$.

\textbf{Necessity:}
If $\widetilde{G}$ is a cycle, by Lemma \ref{le:2.3}, then $\widetilde{G}=\widetilde{C}_{g}$ is of Type $4$.

We assume that $\widetilde{G}$ is not a cycle and $\widetilde{C}_{g}$ is the shortest cycle in $\widetilde{G}$.

\textbf{Case 1.} $g=3$.

In this case, $r(\widetilde{G})=2$. Combining $g=3$ and Theorem \ref{th:4.10}, the reduced graph of $\widetilde{G}$ is $\widetilde{C}_{3}$, which is of Type $4$.

\textbf{Case 2.} $g=4$.

Combining $g=4$, $r(\widetilde{G})=3$ and Lemma \ref{le:2.6},
we have that the reduced graph $\widetilde{G}_{1}$ of $\widetilde{G}$ is a bicyclic graph with rank $3$ and without pendant vertices.
But according to Theorem \ref{th:4.4}, we know that there is no bicyclic graph with $g=4$ and $r(\widetilde{G}_{1})=3$.

\textbf{Case 3.} $g\geq5$.

Since $\widetilde{G}$ is not a cycle, there exists a vertex $x\in V(\widetilde{G})\setminus V(\widetilde{C}_{g})$ such that $d_{\widetilde{C}_{g}}(x)\geq1$.
Since $g\geq5$, by Lemma \ref{le:2.9}, we have $d_{\widetilde{C}_{g}}(x)=1$.
By Lemmas \ref{le:2.1} and \ref{le:2.6},
$$r(\widetilde{C}_{g}+x)=2+r(\widetilde{P}_{g-1})\geq 2+g-2=g>r(\widetilde{G}),$$
a contradiction. \quad $\square$

\begin{figure}[htbp]
  \centering
  \includegraphics[scale=1.0]{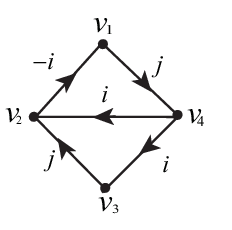}
  \caption{$\widetilde{G}$.}
\end{figure}

\noindent\begin{example}
Consider the $U(\mathbb{Q})$-gain graph $\widetilde{G}$ in Fig. 6, then
$\varphi(v_{1}v_{2}v_{4}v_{1})=i(-i)(-j)=-j$, $\varphi_{v_{1}v_{2}}=-k\varphi_{v_{3}v_{2}}$ and $\varphi_{v_{1}v_{4}}=-k\varphi_{v_{3}v_{4}}$.
Thus $v_{1},v_{3}$ are multiple vertices of $\widetilde{G}$.
Therefore, the subgraph induced $\widetilde{C}_{3}$ on vertices $v_{1},v_{2},v_{4}$ is the reduced graph of $\widetilde{G}$.
By Definition \ref{de:2.2}, $\widetilde{C}_{3}$ is of Type $4$.
Combining Lemmas \ref{le:2.3}, \ref{le:4.8} and $g=3$, we have $r(\widetilde{G})=r(\widetilde{C}_{3})=2=g-1$.
Thus $\widetilde{G}$ is an extremal graph which satisfies the condition in Theorem \ref{th:5.1}.
\end{example}

\section{$U(\mathbb{Q})$-gain graphs $\widetilde{G}$ with $r(\widetilde{G})=g$}

In this section, we will  characterize the graphs that satisfy $r(\widetilde{G})=g$.

\begin{figure}[htbp]
  \centering
  \includegraphics[scale=0.8]{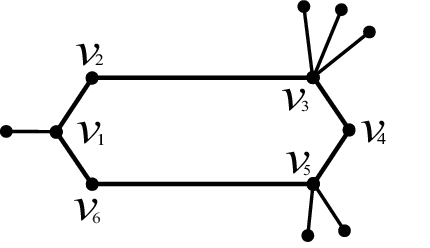}
  \caption{A canonical unicyclic graph $\widetilde{G}$.}
\end{figure}

If $\widetilde{G}$ is a cycle and $m$ pendant stars are attached to its $m$ vertices (where $m$ can be $0$, meaning no vertices are attached), then $\widetilde{G}$ is called a \emph{canonical unicyclic graph} (given in Fig. 7).

\noindent\begin{lemma}\label{le:4.6}
Let $\widetilde{G}$ be a $U(\mathbb{Q})$-gain canonical unicyclic graph of order $n$ with girth $g$, let $t$ be the number of vertices on the cycle of $\widetilde{G}$ that have an attached pendant star, and let $k$ be the number of paths of even (possibly zero) order among the $t$ resulting paths when the pendant stars are all removed. If $\widetilde{G}\neq\widetilde{C}_{n}$, then $t\geq1$, $r(\widetilde{G})=g+k$ and $g\equiv k(\emph{mod}~2)$.
\end{lemma}
\noindent\textbf{Proof.}
By the definition of the canonical unicyclic graph and $\widetilde{G}\neq\widetilde{C}_{n}$, we know that $t\geq1$.
Let $n_{1},n_{2},\ldots,n_{t}$ be the orders of the $t$ pendant stars attached to the cycle of $\widetilde{G}$ and let $m_{1},m_{2},\ldots,m_{t}$ be the orders of the paths obtained upon removing the $t$ pendant stars from
$\widetilde{G}$.
By Lemma \ref{le:2.1}, $\sum_{i=1}^{t}r(\tilde{P}_{m_{i}})=\sum_{i=1}^{t}m_{i}-t+k=g-2t+k$.
By Lemmas \ref{le:2.6} and \ref{le:2.7}, we obtain
$$r(\widetilde{G})=\sum_{i=1}^{t}r(\widetilde{P}_{m_{i}})+2t=g+k.$$
Let $m_{i_{1}},\ldots,m_{i_{k}}$ be even and let $m_{i_{k+1}},\ldots,m_{i_{t}}$ be odd. Then
\begin{align*}
g&=t+m_{1}+\cdots+m_{t}\\
&=k+(t-k)+m_{i_{1}}+\cdots+m_{i_{k}}+m_{i_{k+1}}+\cdots+m_{i_{t}}\\
&=k+m_{i_{1}}+\cdots+m_{i_{k}}+(m_{i_{k+1}}+1)+\cdots+(m_{i_{t}}+1).
\end{align*}
Thus we have $g\equiv k(\textrm{mod}~2)$.\quad $\square$

\noindent\begin{lemma}\label{le:4.7}
Let $\widetilde{G}$ be a $U(\mathbb{Q})$-gain canonical unicyclic graph of order $n$ with girth $g$ and different from a cycle. Then $r(\widetilde{G})=g$ if and only if $g$ is even, there are $t$ ($t\geq1$) vertices on the cycle that have an attached pendant star and the distance between any two such vertices must be even.
\end{lemma}
\noindent\textbf{Proof.}
Let $t$ and $k$ have the same meanings as given in Lemma \ref{le:4.6}.

\textbf{Sufficiency:} In this case, $k=0$. By Lemma \ref{le:4.6}, $r(\widetilde{G})=g$.

\textbf{Necessity:} Combining $r(\widetilde{G})=g$ and Lemma \ref{le:4.6}, we have $k=0$ and $g\equiv 0(\textrm{mod}~2)$. Thus $g$ is even, there are $t$ ($t\geq1$) vertices on the cycle that have an attached pendant star and the distance between any two such vertices must be even.
\quad $\square$\\

Now, we will characterize the $U(\mathbb{Q})$-gain graphs with $r(\widetilde{G})=g$.

\noindent\begin{theorem}\label{th:4.9}
Let $\widetilde{G}$ be a connected $U(\mathbb{Q})$-gain graph with girth $g=3$. Then $r(\widetilde{G})=g$ if and only if $\widetilde{G}$ is one of the following types:
\begin{enumerate}[(a)]
\item $\widetilde{G}=\widetilde{C}_{3}$ is of Type $3$;
\item The reduced graph of $\widetilde{G}$ is $\widetilde{C}_{3}$, which is of Type $3$.
\end{enumerate}
\end{theorem}
\noindent\textbf{Proof.}
\textbf{Sufficiency:}
If $\widetilde{G}=\widetilde{C}_{3}$ is of Type $3$, by Lemma \ref{le:2.3}, then $r(\widetilde{G})=3=g$.

If the reduced graph of $\widetilde{G}$ is $\widetilde{C}_{3}$, which is of Type $3$, then by Lemmas \ref{le:2.3} and \ref{le:4.8}, we have $r(\widetilde{G})=r(\widetilde{C}_{3})=3=g$.

\textbf{Necessity:} If $\widetilde{G}$ is a cycle, by Lemma \ref{le:2.3}, then $\widetilde{G}=\widetilde{C}_{3}$ is of Type $3$.

Assume that $\widetilde{G}$ is not a $U(\mathbb{Q})$-gain cycle and $\widetilde{C}_{3}$ is the shortest cycle in $\widetilde{G}$.
Choosing a vertex $x_{0}\in V(\widetilde{G})$.
Let $Y=N_{\widetilde{G}}(x_{0})$ and $X=V(\widetilde{G})\backslash Y$.
Since $\widetilde{G}$ is connected, $Y$ contains at least one vertex $y_{0}$.
If there exists two vertices $x_{1},x_{2}\in X\backslash\{x_{0}\}$ are adjacent, by Lemma \ref{le:2.6}, then $r(x_{0}+y_{0}+x_{1}+x_{2})=4>r(\widetilde{G})$, a contradiction.
So any vertices in $X$ are not adjacent.
Next, we will prove that $x$ is adjacent to $y$ for any $x\in X$, $y\in Y$.
If there exists $x'\in X$ that is not adjacent to $y'\in Y$, then there exists $y''\in Y$ that is adjacent to $x'$ (since $\widetilde{G}$ is connected).
By Lemma \ref{le:2.6}, $r(x_{0}+x'+y'+y'')=4>r(\widetilde{G})$, a contradiction.

We claim that the induced subgraph $\widetilde{G}_{1}$ with vertex set $Y$ must be a $U(\mathbb{Q})$-gain bipartite graph.
Otherwise, $\widetilde{G}_{1}$ contains a odd cycle $\widetilde{C}_{q}$ as an induced subgraph.
Combining $3=r(\widetilde{G})\geq r(\widetilde{C}_{q})$ and Lemma \ref{le:2.3}, we have $q=3$.
Then $\widetilde{G}$ has an induced subgraph $\widetilde{K}_{4}$.
By Lemma \ref{le:4.11}, $r(\widetilde{K}_{4})=4>r(\widetilde{G})$, a contradiction.
Next, we will prove that $\widetilde{G}_{1}$ is a $U(\mathbb{Q})$-gain complete bipartite graph.
Otherwise, combining $g=3$, we obtain that $\widetilde{G}_{1}$ consists of some $\widetilde{P}_{2}$ and isolated vertices.
Let $\widetilde{P}_{2}$ be a induced subgraph of $\widetilde{G}_{1}$ and let $v\in Y$, $u\in X$.
By Lemma \ref{le:2.6}, $r(\widetilde{P}_{2}+u+v)=4>r(\widetilde{G})$, a contradiction.
Hence, $\widetilde{G}_{1}$ must be a $U(\mathbb{Q})$-gain complete bipartite graph $\widetilde{K}_{s,t}$ with vertex partition ($Y_{1},Y_{2}$), where $|Y_{1}|=s$ and $|Y_{2}|=t$.
Therefore, $\widetilde{G}$ is a $U(\mathbb{Q})$-gain complete tripartite graph $\widetilde{K}_{r,s,t}$ with vertex partition ($X,Y_{1},Y_{2}$), where $|X|=r$.
We may assume that
\begin{align*}
 &A(\widetilde{K}_{r,s,t})=\left (
 \begin{array}{ccccccc}
 \textbf{0} & A_{XY_{1}} & A_{XY_{2}}\\
 A_{XY_{1}}^{\ast} & \textbf{0} & A_{Y_{1}Y_{2}}\\
 A_{XY_{2}}^{\ast} & A_{Y_{1}Y_{2}}^{\ast} & \textbf{0}\\
 \end{array}
 \right)
 =\left (
 \begin{array}{ccccccc}
 A_{1}\\
 A_{2}\\
 A_{3}\\
 \end{array}
 \right),\\
  \end{align*}
where $A_{XY_{1}}$ is a submatrix representing the connection relationship between vertex set $X$ and vertex set $Y_{1}$, and
$A_{XY_{1}}^{\ast}$ (resp., $A_{XY_{2}}^{\ast}$, $A_{Y_{1}Y_{2}}^{\ast}$) is the conjugate transpose of $A_{XY_{1}}$ (resp., $A_{XY_{2}}$, $A_{Y_{1}Y_{2}}$).
We know that $A_{XY_{1}}$, $A_{XY_{2}}$ and $A_{Y_{1}Y_{2}}$ contain nonzero quaternion as entries.
If there exists a row submatrix $A_{i}$ ($i\in\{1,2,3\}$) such that $r(A_{i})\geq 2$, then $r(\widetilde{K}_{r,s,t})\geq 2r(A_{i})\geq 4$, a contradiction.
So for any $1\leq i\leq3$, each row vector of $A_{i}$ is the nonzero multiple of the first row vector.
It follows that the reduced graph of $\widetilde{K}_{r,s,t}$ is $\widetilde{C}_{3}$.
Combining $r(\widetilde{G})=3$ and Lemma \ref{le:4.8}, $r(\widetilde{C}_{3})=r(\widetilde{G})=3$.
By Lemma \ref{le:2.3}, $\widetilde{C}_{3}$ is of Type $3$.
\quad $\square$

\noindent\begin{theorem}\label{th:5.10}
Let $\widetilde{G}$ be a connected $U(\mathbb{Q})$-gain graph with girth $g=4$. Then $r(\widetilde{G})=g$ if $\widetilde{G}$ is one of the following types:
\begin{enumerate}[(a)]
\item $\widetilde{G}=\widetilde{C}_{4}$ is of Type $2$;
\item $\widetilde{G}$ is $\widetilde{C}_{4}$ with a vertex or two vertices (that are not adjacent) attached as a pendant star;
\item $\widetilde{G}$ is obtained from the cycle $\widetilde{C}_{4}$ which is of Type $1$ and the star $\widetilde{S}_{q+1}$ ($q\geq1$) by joining a vertex of the cycle to the center of the star;
\item The reduced graph of $\widetilde{G}$ is switching equivalent to one of the $U(\mathbb{Q})$-gain graphs with rank $4$ listed in Table $1$ except $\widetilde{G}_{1}$, $\widetilde{G}_{2}$, $\widetilde{G}_{5}$, $\widetilde{G}_{6}$ and $\widetilde{G}_{7}$;
\item The reduced graph of $\widetilde{G}$ is switching equivalent to one of the $\widetilde{G}_{18}-\widetilde{G}_{22}$ listed in Table $2$.
\end{enumerate}
\end{theorem}
\noindent\textbf{Proof.}
If $\widetilde{G}=\widetilde{C}_{4}$ is of Type $2$, by Lemma \ref{le:2.3}, then $r(\widetilde{G})=4=g$.

If  $\widetilde{G}$ is $\widetilde{C}_{4}$ with a vertex or two vertices (that are not adjacent) attached as a pendant star, by Lemma \ref{le:2.6}, then $r(\widetilde{G})=g$.

If $\widetilde{G}$ is obtained from the cycle $\widetilde{C}_{4}$ which is of Type $1$ and the star $\widetilde{S}_{q+1}$ ($q\geq1$) by joining a vertex of the cycle to the center of the star, by Lemmas \ref{le:2.3} and \ref{le:2.6}, then $r(\widetilde{G})=r(\widetilde{C}_{4})+2=g$.

If the reduced graph of $\widetilde{G}$ is switching equivalent to one of the $U(\mathbb{Q})$-gain graphs with rank $4$ listed in Table $1$ except $\widetilde{G}_{1}$, $\widetilde{G}_{2}$, $\widetilde{G}_{5}$, $\widetilde{G}_{6}$ and $\widetilde{G}_{7}$,
or one of the $\widetilde{G}_{18}-\widetilde{G}_{22}$ listed in Table $2$, by Theorems \ref{th:4.4}, \ref{th:4.5} and Lemma \ref{le:4.8}, then $r(\widetilde{G})=4=g$.
\quad $\square$

\begin{figure}[htbp]
  \centering
  \includegraphics[scale=0.75]{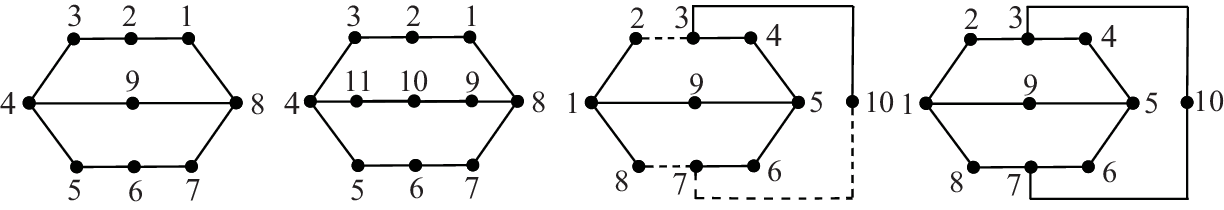}
  \caption{$\widetilde{\theta}(1,3,3)$, $\widetilde{\theta}(3,3,3)$, $\widetilde{T}_{0}$ and $\widetilde{G}_{0}$.}
\end{figure}

Next, we will characterize three $U(\mathbb{Q})$-gain graphs with ranks equal to their girths.

\noindent\begin{lemma}\label{le:5.12}
\begin{enumerate}[(a)]
\item $r(\widetilde{\theta}(1,3,3))=6$ if and only if all the cycles in $\widetilde{\theta}(1,3,3)$ are of Type $1$;
\item $r(\widetilde{\theta}(3,3,3))=8$ if and only if all the cycles in $\widetilde{\theta}(3,3,3)$ are of Type $1$;
\item $r(\widetilde{G}_{0})=6$ (see Fig. 8) if and only if all the cycles in $\widetilde{G}_{0}$ are of Type $1$.
\end{enumerate}
\end{lemma}
\noindent\textbf{Proof.}
If $\widetilde{G}=\widetilde{\theta}(1,3,3)$, then the switching function $\xi_{2}$ and the spanning tree in $\widetilde{G}^{\xi_{2}}$ is similar to those in $\widetilde{G}_{b}^{\xi_{b}}$ in Fig. 3.
Thus $\varphi^{\xi_{2}}(v_{4}v_{9})=\varphi(v_{1}v_{2}v_{3}v_{4}v_{9}v_{8}v_{1})$ and $\varphi^{\xi_{2}}(v_{7}v_{8})=\varphi(v_{1}v_{2}v_{3}v_{4}v_{5}v_{6}v_{7}v_{8}v_{1})$.
Let $a=\varphi^{\xi_{2}}(v_{4}v_{9})$ and $b=\varphi^{\xi_{2}}(v_{7}v_{8})$.
After calculation, we get
\begin{align*}
 r(\widetilde{G}^{\xi_{2}})=6+\left
 (\begin{array}{cccccccc}
 0 & b-1 & 0\\
 \overline{b}-1 & 0 & a+1\\
 0 & \overline{a}+1 & 0\\
\end{array}
 \right).
 \end{align*}
So $r(\widetilde{G}^{\xi_{2}})=6$ if and only if $a=-1$ and $b=1$, i.e., $\varphi(v_{1}v_{2}v_{3}v_{4}v_{9}v_{8}v_{1})=-1$ and $\varphi(v_{1}v_{2}v_{3}v_{4}$ $v_{5}v_{6}v_{7}v_{8}v_{1})=1$. Thus $\varphi(v_{4}v_{5}v_{6}v_{7}v_{8}v_{9}v_{4})=-1$.
By Definition \ref{de:2.2}, all the cycles in $\widetilde{G}^{\xi_{2}}$ are of Type $1$.
Therefore, $r(\widetilde{\theta}(1,3,3))=6$ if and only if all the cycles in $\widetilde{\theta}(1,3,3)$ are of Type $1$.

If $\widetilde{G}=\widetilde{\theta}(3,3,3)$, then the switching function $\xi_{2}$ and the spanning tree in $\widetilde{G}^{\xi_{2}}$ is similar to those in $\widetilde{G}_{b}^{\xi_{b}}$ in Fig. 3.
Thus $\varphi^{\xi_{2}}(v_{4}v_{11})=\varphi(v_{1}v_{2}v_{3}v_{4}v_{11}v_{10}v_{9}v_{8}v_{1})$ and $\varphi^{\xi_{2}}(v_{7}v_{8})=\varphi(v_{1}v_{2}v_{3}v_{4}v_{5}v_{6}v_{7}v_{8}v_{1})$.
Let $a=\varphi^{\xi_{2}}(v_{4}v_{11})$ and $b=\varphi^{\xi_{2}}(v_{7}v_{8})$.
After calculation, we get
\begin{align*}
 r(\widetilde{G}^{\xi_{2}})=8+\left
 (\begin{array}{cccccccc}
 0 & b-1 & 0\\
 \overline{b}-1 & 0 & a-1\\
 0 & \overline{a}-1 & 0\\
\end{array}
 \right).
 \end{align*}
So $r(\widetilde{G}^{\xi_{2}})=8$ if and only if $a=b=1$, i.e., $\varphi(v_{1}v_{2}v_{3}v_{4}v_{11}v_{10}v_{9}v_{8}v_{1})=\varphi(v_{1}v_{2}v_{3}v_{4}v_{5}v_{6}v_{7}v_{8}v_{1})=1$. Thus $\varphi(v_{4}v_{5}v_{6}v_{7}v_{8}v_{9}v_{10}v_{11}v_{4})=1$.
By Definition \ref{de:2.2}, all the cycles in $\widetilde{G}^{\xi_{2}}$ are of Type $1$.
Therefore, $r(\widetilde{\theta}(3,3,3))=8$ if and only if all the cycles in $\widetilde{\theta}(3,3,3)$ are of Type $1$.

If $\widetilde{G}=\widetilde{G}_{0}$, then let $\widetilde{T}_{0}$ (see Fig. 8) be a spanning tree of $\widetilde{G}_{0}$ and let $v_{1}$ be a root vertex.
For $v_{i}\in V(\widetilde{G}_{0})\backslash\{v_{1}\}$ ($i=2,3,\ldots,10$), let $\widetilde{P}_{v_{i}v_{1}}$ be the unique path from $v_{i}$ to $v_{1}$ in $\widetilde{T}_{0}$.
Let $\xi_{3}(v_{i})=\varphi(\widetilde{P}_{v_{i}v_{1}})$ and $\xi_{3}(v_{1})=1$.
After using the switching function $\xi_{3}$, we have that all edges in $\widetilde{T}_{0}$ have gain $1$.
So in $\widetilde{G}_{0}^{\xi_{3}}$,
$\varphi^{\xi_{3}}(v_{2}v_{3})=\varphi(v_{1}v_{2}v_{3}v_{4}v_{5}v_{9}v_{1})$,
$\varphi^{\xi_{3}}(v_{7}v_{8})=\varphi(v_{1}v_{9}v_{5}v_{6}v_{7}v_{8}v_{1})$, $\varphi^{\xi_{3}}(v_{7}v_{10})=\varphi(v_{3}v_{4}v_{5}v_{6}v_{7}v_{10}v_{3})$
and the gains of the other edges have gain $1$.
Let $a=\varphi^{\xi_{3}}(v_{2}v_{3})$, $b=\varphi^{\xi_{3}}(v_{7}v_{8})$ and $c=\varphi^{\xi_{3}}(v_{7}v_{10})$.
After calculation, we get
\begin{align*}
 r(\widetilde{G}_{0}^{\xi_{3}})=6+\left
 (\begin{array}{cccccccc}
 0 & -\overline{a}+b & -\overline{a}-1 & c+1\\
 \overline{b}-a & 0 & 0 & 0 \\
\end{array}
 \right).
 \end{align*}
So $r(\widetilde{G}_{0}^{\xi_{3}})=6$ if and only if $a=b=c=-1$, i.e., all the cycles in $\widetilde{G}_{0}^{\xi_{3}}$ are of Type $1$.
\quad $\square$

\noindent\begin{theorem}\label{th:5.11}
Let $\widetilde{G}$ be a connected $U(\mathbb{Q})$-gain graph with girth $g\geq5$. Then $r(\widetilde{G})=g$ if and only if $\widetilde{G}$ is one of the following types:
\begin{enumerate}[(a)]
\item $\widetilde{G}=\widetilde{C}_{g}$ is of Type $2$ or $3$;
\item $\widetilde{G}=\widetilde{\theta}(1,3,3)$ and all the cycles in $\widetilde{\theta}(1,3,3)$ are of Type $1$;
\item $\widetilde{G}=\widetilde{\theta}(3,3,3)$ and all the cycles in $\widetilde{\theta}(3,3,3)$ are of Type $1$;
\item $\widetilde{G}=\widetilde{G}_{0}$ (see Fig. 8) and all the cycles in $\widetilde{G}_{0}$ are of Type $1$;
\item A $U(\mathbb{Q})$-gain canonical unicyclic graph with an even girth $g$, there are $t$ ($1\leq t\leq g$) vertices on the cycle that have an attached pendant star and the distance between any two such vertices must be even;
\item A $U(\mathbb{Q})$-gain unicyclic graph obtained from the cycle $\widetilde{C}_{g}$ which is of Type $1$ and the star $\widetilde{S}_{q+1}$ ($q\geq1$) by joining a vertex of the cycle to the center of the star.
\end{enumerate}
\end{theorem}
\noindent\textbf{Proof.}
\textbf{Sufficiency:}
If $\widetilde{G}=\widetilde{C}_{g}$ is of Type $2$ or $3$, by Lemma \ref{le:2.3}, then $r(\widetilde{G})=g$.

If $\widetilde{G}=\widetilde{\theta}(1,3,3)$ and all the cycles in $\widetilde{\theta}(1,3,3)$ are of Type $1$,
or $\widetilde{G}=\widetilde{\theta}(3,3,3)$ and all the cycles in $\widetilde{\theta}(3,3,3)$ are of Type $1$,
or $\widetilde{G}=\widetilde{G}_{0}$ (see Fig. 8) and all the cycles in $\widetilde{G}_{0}$ are of Type $1$,
by Lemma \ref{le:5.12}, then $r(\widetilde{G})=g$.

If $\widetilde{G}$ is a $U(\mathbb{Q})$-gain canonical unicyclic graph with an even girth $g$, there are $t$ ($t\geq1$) vertices on the cycle that have an attached pendant star and the distance between any two such vertices must be even, by Lemma \ref{le:4.7}, then $r(\widetilde{G})=g$.

If $\widetilde{G}$ is a $U(\mathbb{Q})$-gain unicyclic graph obtained from the cycle $\widetilde{C}_{g}$ which is of Type $1$ and the star $\widetilde{S}_{q+1}$ ($q\geq1$) by joining a vertex of the cycle to the center of the star, by Lemmas \ref{le:2.3} and \ref{le:2.6}, then $r(\widetilde{G})=r(\widetilde{C}_{g})+2=g$.

\textbf{Necessity:} If $\widetilde{G}$ is a cycle, by Lemma \ref{le:2.3}, then $\widetilde{G}=\widetilde{C}_{g}$ is of Type $2$ or $3$.

Assume that $\widetilde{G}$ is not a $U(\mathbb{Q})$-gain cycle and $\widetilde{C}_{g}$ is the shortest cycle in $\widetilde{G}$.
Let $Y_{t}=\{y\in V(\widetilde{G})\backslash V(\widetilde{C}_{g})~|~d_{\widetilde{G}}(y,\widetilde{C}_{g})=t\}$.
Since $g\geq5$ and Lemma \ref{le:2.9}, $|N_{\widetilde{C}_{g}}(v)|=1$ for every $v\in Y_{1}$.
If $Y_{3}\neq\emptyset$, then we find a vertex $y\in Y_{3}$. Let $\widetilde{P}$ be the shortest path from $y$ to $\widetilde{C}_{g}$. By Lemmas \ref{le:2.6} and \ref{le:2.1},
$$r(\widetilde{C}_{g}+\widetilde{P})=4+r(\widetilde{P}_{g-1})\geq g+2>r(\widetilde{G}),$$
a contradiction. Thus $Y_{i}=\emptyset~(i\geq3)$.

\textbf{Case 1.} $|Y_{2}|=0$.

In this case, $\widetilde{G}$ is a canonical unicyclic graph.
By Lemma \ref{le:4.7}, $g$ is even, there are $t$ ($1\leq t\leq g$) vertices on the cycle that have an attached pendant star and the distance between any two such vertices must be even.

\textbf{Case 2.} $|Y_{2}|=1$.

Let $Y_{2}=\{x\}$. If $|Y_{1}|=1$, by Lemmas \ref{le:2.3} and \ref{le:2.6}, then $\widetilde{G}$ is a $U(\mathbb{Q})$-gain unicyclic graph obtained from the cycle $\widetilde{C}_{g}$ which is of Type $1$ and $\widetilde{P}_{2}$ by joining a vertex of the cycle to a vertex of $\widetilde{P}_{2}$.
So we assume that $|Y_{1}|\geq2$.
Let $N_{Y_{1}}(x)=\{y\}$ and $y'\in Y_{1}~(y'\neq y)$. If $x$ is not adjacent to $y'$, by Lemmas \ref{le:2.1} and \ref{le:2.6}, then
$$r(\widetilde{C}_{g}+x+y+y')=4+r(\widetilde{P}_{g-1})\geq g+2>r(\widetilde{G}),$$
a contradiction. Thus $x$ is adjacent to $y'$. Let $V(\widetilde{C}_{g})=\{v_{1},v_{2},\ldots,v_{g}\}$, $N_{\widetilde{C}_{g}}(y)=\{v_{1}\}$ and $N_{\widetilde{C}_{g}}(y')=\{v_{i}\}$.
If $g$ is odd, by Lemmas \ref{le:2.1} and \ref{le:2.6}, then $r(\widetilde{C}_{g}+y)=2+r(\widetilde{P}_{g-1})= g+1>r(\widetilde{G})$, a contradiction. Thus $g$ is even.
Without loss of generality, $i\leq\frac{g}{2}+1$.
Then $\frac{g}{2}+4\geq(i-1)+4\geq g$, thus $g=6$ and $i=3$, or $g=6$ and $i=4$, or $g=8$ and $i=5$.

\textbf{Case 2.1.} $g=6$ and $i=3$.

In this case, $\widetilde{G}=\widetilde{\theta}(1,3,3)$ or $\widetilde{G}=\widetilde{G}_{0}$ (see Fig. 8).
If $\widetilde{G}=\widetilde{\theta}(1,3,3)$, by Lemma \ref{le:5.12}, then each cycle in $\widetilde{\theta}(1,3,3)$ is of Type $1$.
If $\widetilde{G}=\widetilde{G}_{0}$, by Lemma \ref{le:5.12}, then each cycle in $\widetilde{G}_{0}$ is of Type $1$.

\textbf{Case 2.2.} $g=6$ and $i=4$.

In this case, $\widetilde{G}=\widetilde{\theta}(2,2,3)$. By Lemmas \ref{le:2.4} and \ref{le:2.6}, $r(\widetilde{G})\geq8>g$, a contradiction.

\textbf{Case 2.3.} $g=8$ and $i=4$.

In this case, $\widetilde{G}=\widetilde{\theta}(3,3,3)$.
By Lemma \ref{le:5.12}, each cycle in $\widetilde{\theta}(3,3,3)$ is of Type $1$.

\textbf{Case 3.} $|Y_{2}|\geq2$.

We first prove that for any distinct vertices $u,v\in Y_{2}$, $N_{Y_{1}}(u)=N_{Y_{1}}(v)$. On the contrary, we find two distinct vertices $u_{1}, v_{1}\in Y_{1}$ such that $u_{1}\in N_{Y_{1}}(u)$ and $v_{1}\in N_{Y_{1}}(v)$.

If $u$ is adjacent to $v$, then $u$ is not adjacent to $v_{1}$ (since $g\geq5$). By Lemmas \ref{le:2.1} and \ref{le:2.6},
$$r(\widetilde{C}_{g}+u+v+v_{1})=4+r(\widetilde{P}_{g-1})\geq g+2>r(\widetilde{G}),$$
a contradiction.

If $u$ is not adjacent to $v$, then exactly one or none of $u$ is adjacent to $v_{1}$ and $v$ is adjacent to $u_{1}$ is satisfied. By Lemmas \ref{le:2.6} and \ref{le:2.3},
$$r(\widetilde{C}_{g}+u+v+u_{1}+v_{1})=4+r(\widetilde{C}_{g})\geq 4+(g-2)=g+2>r(\widetilde{G}),$$
a contradiction.

Thus exists $y\in Y_{1}$ such that $N_{Y_{1}}(x)=\{y\}$ for all $x\in Y_{2}$. Next, we prove that $|Y_{1}|=1$.
On the contrary, there exists $y'\in Y_{1}$ such that $y\neq y'$. Choosing any $x\in Y_{2}$, then $x$ is adjacent to $y$ and $x$ is not adjacent to $y'$. By Lemmas \ref{le:2.1} and \ref{le:2.6},
$$r(\widetilde{C}_{g}+x+y+y')=4+r(\widetilde{P}_{g-1})\geq g+2>r(\widetilde{G}),$$
a contradiction. Thus $|Y_{1}|=1$. By Lemma \ref{le:2.6}, $r(\widetilde{G})=2+r(\widetilde{C}_{g})=g$, then $r(\widetilde{C}_{g})=g-2$.
By Lemma \ref{le:2.3}, $\widetilde{C}_{g}$ is of Type $1$. Let $|Y_{2}|=q\geq2$, then $\widetilde{G}$ is a $U(\mathbb{Q})$-gain unicyclic graph obtained from the cycle $\widetilde{C}_{g}$ which is of Type $1$ and the star $\widetilde{S}_{q+1}$ ($q\geq2$) by joining a vertex of the cycle to the center of the star.
\quad $\square$

\begin{figure}[htbp]
  \centering
  \includegraphics[scale=1.0]{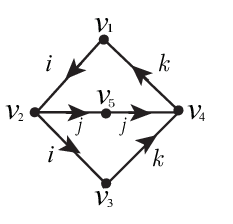}
  \caption{$\widetilde{G}$.}
\end{figure}

\noindent\begin{example}
Consider the $U(\mathbb{Q})$-gain graph $\widetilde{G}=\widetilde{\theta}(1,1,1)$ in Fig. 9, then
$$\varphi(v_{1}v_{2}v_{5}v_{4}v_{1})=ijjk=kjk=-ik=j;$$
$$\varphi(v_{1}v_{2}v_{3}v_{4}v_{1})=iikk=-kk=1.$$
By Definition \ref{de:2.2}, the subgraph induced on vertices $v_{1},v_{2},v_{5},v_{4}$ is of Type $2$ and the subgraph induced on vertices $v_{1},v_{2},v_{3},v_{4}$ is of Type $1$.
Combining Theorem \ref{th:4.4} and $g=4$, $r(\widetilde{G})=4=g$.
Thus $\widetilde{G}$ is an extremal graph which satisfies the condition in Theorem \ref{th:5.10}.
\end{example}

\noindent\begin{remark}
In Theorem \ref{th:5.10}, the conditions of this theorem are sufficient but not necessary. The problem is that all quaternion unit  gain graphs with rank $4$ have not yet been characterized. This is a problem we will continue to consider later.
\end{remark}

\textbf{Declaration of competing interest}

The authors declare that they have no known competing financial interests or personal relationships that could have appeared to influence the work reported in this paper.

\textbf{Date availability}

No date was used for the research described in the article.

\end{document}